\documentclass[11pt]{article}

\usepackage{amsmath}
\usepackage{amssymb,bbm}
\usepackage{mathabx}
 \usepackage{IEEEtrantools}
\usepackage{mathtools}

\usepackage{mathrsfs}
\usepackage{url}
\usepackage{hyperref}
\usepackage{xcolor}
\usepackage{color}

\usepackage{graphicx}

\usepackage{multirow}

\usepackage{authblk}

\usepackage{pgf,tikz}
\usepackage{tikz-3dplot}
\usetikzlibrary{arrows}
\usetikzlibrary{arrows.meta}

\newcommand{\sgn}{\mbox{sign\,}}

\newcommand{\prf}{\noindent{\it Proof\/}: }

\def \qed { \mbox{}\hfill
$\Box$\vspace{1ex}}

\newcommand{\half}{\frac{1}{2}}
\newcommand{\halfof}[1]{\frac{#1}{2}}

\counterwithin*{equation}{section}

\newtheorem{dfn}{Definition}[section]
\newtheorem{thm}[dfn]{Theorem}
\newtheorem{lmma}[dfn]{Lemma}
\newtheorem{ppsn}[dfn]{Proposition}
\newtheorem{crlre}[dfn]{Corollary}
\newtheorem{xmpl}[dfn]{Example}
\newtheorem{rmrk}[dfn]{Remark}
\newtheorem{xrcs}[dfn]{Exercise}

\newcommand{\bdfn}{\begin{dfn}\rm}
\newcommand{\bthm}{\begin{thm}}
\newcommand{\blmma}{\begin{lmma}}
\newcommand{\bppsn}{\begin{ppsn}}
\newcommand{\bcrlre}{\begin{crlre}}
\newcommand{\bxmpl}{\begin{xmpl}}
\newcommand{\brmrk}{\begin{rmrk}\rm}
\newcommand{\bxrcs}{\begin{xrcs}\rm\footnotesize}

\newcommand{\edfn}{\end{dfn}}
\newcommand{\ethm}{\end{thm}}
\newcommand{\elmma}{\end{lmma}}
\newcommand{\eppsn}{\end{ppsn}}
\newcommand{\ecrlre}{\end{crlre}}
\newcommand{\exmpl}{\end{xmpl}}
\newcommand{\ermrk}{\end{rmrk}}
\newcommand{\exrcs}{\end{xrcs}}

\begin{document}
%%--------------------------------------------------------------------
\author[1]{\sc Partha Sarathi Chakraborty\thanks{parthacsarathi.isi.smu@gmail.com, parthac@imsc.res.in}}
\author[2]{\sc Arup Kumar Pal\thanks{Supported by SERB MATRICS grant MTR/2017/000544}\thanks{arup@isid.ac.in, arupkpal@gmail.com}}
\affil[1]{Indian Statistical Institute,  Kolkata, INDIA}
\affil[2]{Indian Statistical Institute, Delhi, INDIA}
\title{An approximate equivalence
for the GNS representation of the Haar state of $SU_{q}(2)$}
\maketitle
%%--------------------------------------------------------------------
%%--------------------------------------------------------------------
%  ABSTRACT
%%--------------------------------------------------------------------
 \begin{abstract} 
 We use the crystallised $C^*$-algebra $C(SU_{q}(2))$ at $q=0$ to obtain a unitary that gives an approximate equivalence involving the GNS
  representation on the $L^{2}$ space of the Haar state of the quantum
  $SU(2)$ group and the direct integral of all the infinite dimensional
  irreducible representations of the $C^{*}$-algebra $C(SU_{q}(2))$ for
  nonzero values of the parameter $q$. This approximate equivalence
  gives a $KK$ class  via the Cuntz picture in terms of
  quasihomomorphisms as well as a Fredholm representation of the dual
  quantum group $\widehat{SU_q(2)}$ with coefficients in a
  $C^*$-algebra in the sense of Mishchenko. 
\end{abstract}
\textbf{2020 AMS Subject Classification No.:}
%17B37, %Quantum groups (quantized enveloping algebras) and related deformations
20G42, %Quantum groups (quantized function algebras) and their representations
58B32, % Geometry of quantum groups
46L67, % Quantum groups (operator algebraic aspects)
19K35 % Kasparov theory ($KK$-theory)
\section{Introduction}
%%--------------------------------------------------------------------
The quantum group $SU_{q}(2)$ was studied in the context of Connes' set up of
spectral triples independently by Chakraborty \& Pal \cite{ChaPal-2003aa} and
Dabrowski et al \cite{DabLanSit-2005aa}. Subsequently Connes \cite
{Con-2004ab} studied the equivariant spectral triple in \cite
{ChaPal-2003aa} for $SU_{q}(2)$ in great detail, in particular determined its
dimension spectrum, established regularity and showed that the local index
theorem proved by him and Moscovici applies to this spectral triple. His idea
was then carried forward by a host of authors for studying equivariant specral
triples for cartain quantum groups and their homogeneous spaces that were known
at that time (\cite{SuiDabLan-2005aa}, \cite{DAnDab-2006aa}, \cite
{ChaPal-2008aa}, \cite{PalSun-2010aa}). A few years later, Neshveyev \& Tuset
\cite{NesTus-2010ab} constructed an equivariant Dirac operator for a large
class of quantum groups, namely the $q$-deformations of all simple simply
connected compact Lie groups. 

However, despite the work of Neshveyev \& Tuset setting the stage, further work on
local index formulas did not take off. In particular regularity, discreteness of
dimension spectrum and computation of local cyclic cocycles still remain open.
This calls for a closer look at the papers by Connes \cite{Con-2004ab} and  
Neshveyev \& Tuset \cite{NesTus-2010ab}. The most important idea used in
Connes' paper is the observation that if one substitutes $q=0$ in the
expressions for the actions of the generating elements $\alpha$ and $\beta$ on
$L^{2}(SU_{q}(2))$, one gets a representation of the $C^*$-algebra $C(SU_{q}
(2))$ at $q=0$ (which is just the $C^*$-algebra if one replaces $q=0$ in the
defining relations of $C(SU_{q}(2))$), and a lot of simplifications occur. Connes
was then able to connect the GNS representation $\lambda_{q}$ on the $L^
{2}$ space of the Haar state with another faithful representation  $\pi_{q}$ of
$C(SU_{q}(2)$ that one gets by combining all the irreducible representations
for these $C^*$-algebras, which have been characterized and listed by results
of Soibelman \cite{KorSoi-1998ab}. Following \cite{MatYun-2022aa}, henceforth
we will refer to this representation $\pi_{q}$ as the Soibelman representation.
It is computationally much more tractable and acts on the Hilbert space $\ell^
{2}(\mathbb{Z}\times\mathbb{N})$. The Dirac operator of Neshveyev \& Tuset, on
the other hand, comes from the classical Dirac through a twisting procedure,
and while it has a sound conceptual origin, as Hilbert space operators, both
this Dirac $D$ as well as the $C(G_{q})$ elements are rather difficult to work
with. Therefore, following Connes' idea, it is reasonable to try to make use of
the  Soibelman representation.

With this in mind, we tried to see what is going on behind the computations in
Connes' paper at the $C^*$-algebra level. In the present paper, we make our
observation precise and explicit. In particular, we prove a certain approximate
equivalence between the GNS representation for the quantum $SU(2)$ group and an
ampliation of the Soibelman representation. The key step in proving this
approximate equivalence is obtaining the unitary that gives the equivalence. We
obtain this unitary by using the crystallised quantised function algebra  $C
(SU_{q}(2))$ at $q=0$. We obtain the unitary in  Section~2 and prove the
approximate equivalence in Section~3. In Section~4, we list a few consequences
of this equivalence. The approximate equivalence also gives an element of the
$KK$ group via the Cuntz picture. Further, we extend the notion of a Fredholm
representation of a discrete group introduced by Mishchenko \cite
{Mis-1975pz} and show that the approximate equivalence presents us with an
example of a Fredholm representation for the dual quantum group $\widehat{SU_
{q}(2)}$. Finally we describe how the approximate equivalence sets up a
relation between the equivariant spectral triple constructed in \cite
{ChaPal-2003aa} and another spectral triple studied by the authors in \cite
{ChaPal-2003ab}, which is precisely what lies behind Connes' computations
in \cite{Con-2004aa}.

\textbf{Notation}: 
$\mathcal{H}$ will denote a complex separable Hilbert space.
 $\mathcal{L}(\mathcal{H})$ and $\mathcal{L}(X)$ will denote the space of
 bounded linear maps on a Hilbert space $\mathcal{H}$ and the space of
 adjointable operators on a Hilbert $C^{*}$-module $X$ respectively. Similarly,
 $\mathcal{K}(\mathcal{H})$ and $\mathcal{K}(X)$ will denote the spaces of
 compact operators on them. We will denote the Toeplitz algebra by $\mathscr
 {T}$. For a real number $t$, $t_{+}$ and $t_{-}$ will denote its positive and
 negative parts respectively. Throughout the paper, $q$ will denote a real
 number in the interval $(-1,1)$ that will be assumed to be nonzero unless explicitly stated
 otherwise. For a positive integer $n$, we will denote by $g(n)$ the number
 $\left(1-q^{2n}\right)^{\half}$.

%%--------------------------------------------------------------------
\section{Unitary equivalence at $q=0$}
%%--------------------------------------------------------------------
%%--------------------------------------------------------------------
\subsection{Quantum $SU(2)$ group}
%%--------------------------------------------------------------------
To fix notation, let us give here a very brief description of the quantum
$SU(2)$ group. The $C^{*}$-algebra associated with $SU_q(2)$, usually denoted by
$C(SU_{q}(2))$, is the $C^*$-algebra generated by two elements $\alpha$ and
$\beta$ satisfying the following relations:
%%-------------------------------------
\begin{IEEEeqnarray*}{rClrCl}
\alpha^*\alpha+\beta^*\beta &=& I,\qquad & \alpha\beta &=& q\beta\alpha,\\
\alpha\alpha^*
+q^2\beta\beta^* &=& I,&\alpha\beta^* &=& q\beta^*\alpha,\yesnumber\label{eq:suq2-relations}\\
&&&\beta^*\beta &=& \beta\beta^*.
\end{IEEEeqnarray*}
%%-------------------------------------
We will also use the symbol $A_{q}$ to denote this $C^{*}$-algebra and use
$\alpha_{q}$ and $\beta_{q}$ for the generating elements of $A_{q}$ instead of
just $\alpha$ and $\beta$. The quantum group structure is given by the coproduct
$\Delta$ which is a unital *-homomorphism from $A_{q}$ to $A_{q}\otimes A_{q}$
given by
%%-------------------------------------
\begin{IEEEeqnarray}{rCl}
\Delta(\alpha_{q})&=&\alpha_{q}\otimes\alpha_{q}-
    q\beta_{q}^*\otimes\beta_{q},\label{eq:qcomult-1}\\
\Delta(\beta_{q})&=&\beta_{q}\otimes\alpha_{q}+
    \alpha_{q}^*\otimes\beta_{q}.  \label{eq:qcomult-2}
\end{IEEEeqnarray}
%%-------------------------------------
For two continuous linear functionals $\rho_1$ and $\rho_2$ on $A_{q}$, one defines
their convolution product by: $\rho_1\ast\rho_2(a)=(\rho_1\otimes\rho_2)\Delta(a)$.
It is known \cite{Wor-1987aa} that $A_{q}$ admits a faithful state $h$, called
the Haar state, that satisfies
%%-------------------------------------
\[
h\ast \rho (a) = h(a)\rho(I) = \rho\ast h(a)
\]
%%-------------------------------------
for all continuous linear functionals $\rho$ and all $a\in A_{q}$.
We will be concerned with the GNS space associated with this state.

%%--------------------------------------------------------------------
\subsection{Representation on $\ell^{2}(\mathbb{N})\otimes \ell^{2}(\mathbb{Z})$}
%%--------------------------------------------------------------------
Let $\{e_{n}\}_{n\in\mathbb{N}}$ be the canonical orthonormal basis for
$\ell^{2}(\mathbb{N})$, where $\mathbb{N}$ denotes the set of nonnegative
integers. Let us denote 
by $P_{0}$ the projection onto $\mathbb{C}e_{0}$, 
by $N$ the number operator on
$\ell^{2}(\mathbb{N})$, and by $S$ the left shift:
%%-------------------------------------
\[
Ne_{k}=ke_{k},\qquad Se_{k}=\begin{cases} e_{k-1}
   &\text{if }k\geq 1,\\ 0&\text{if }k=0.\end{cases}
\]
%%-------------------------------------
It is well-known and easy to show that as $z$ varies over the unit circle $S^{1}$,
the following constitute all inequivalent infinite dimensional irreducible
representations of the $C^*$-algebra $A_{q}$:
%%-------------------------------------
\begin{IEEEeqnarray*}{rCl}
	\alpha_{q} &\mapsto& S\sqrt{1-q^{2N}}\\
	\beta_{q} &\mapsto& zq^{N}.
\end{IEEEeqnarray*}
%%-------------------------------------
The direct integral of these representations gives a faithful representation
$\pi_{q}$ of $A_{q}$ on the Hilbert space
$\mathcal{H}_{\pi}:=\ell^{2}(\mathbb{N})\otimes \ell^{2}(\mathbb{Z})$ given by
%%-------------------------------------
\begin{equation}\label{pi}
\pi_{q}(\alpha_{q})=S\sqrt{I-q^{2N}}\otimes I,
    \qquad \pi_{q}(\beta_{q})=q^N\otimes S,
\end{equation}
%%-------------------------------------
where we have used the same symbol $S$ to denote the left shift $e_k\mapsto
e_{k-1}$ on $\ell^{2}(\mathbb{Z})$.

%%--------------------------------------------------------------------
\subsection{GNS representation on the $L^{2}$ space}
%%--------------------------------------------------------------------
Let $\mathcal{H}\equiv L^{2}(SU_q(2))$ denote the GNS space for the Haar state on
$A_{q}$. Let $\lambda_{q}$ denote the GNS representation (i.e.\ the representation
by left multiplication) of $A_{q}$ on $\mathcal{H}$. By the Peter-Weyl theorem for
compact quantum groups, and the representation theory for the quantum group
$SU_q(2)$, it follows that $\mathcal{H}$ has a natural orthonomal basis
$\{e^{n}_{ij}: n\in\half\mathbb{N}, i,j=-n,-n+1,\ldots,n-1,n\}$ where
$e^{n}_{ij}$'s are the normalized matrix entries of the irreducible unitary
(co-)representations of $SU_q(2)$. Thus $\mathcal{H}$ can be identified with
$\ell^2(\Gamma)$ where
%%-------------------------------------
\[
\Gamma=\left\{(n,i,j): n\in  \half\mathbb{N},  \; i,j\in \{-n, -n+1,\ldots,n-1,n\}\right\}.
\]
%%-------------------------------------
The representation $\lambda_q$ of $A_q$ can be written down
explicitly using Clebsch-Gordan coefficients (see Equations (2.1--2.2),
\cite{ChaPal-2003aa}) as follows:
%%-------------------------------------
\begin{IEEEeqnarray}{rCrCl}
\lambda_{q}(\alpha_{q}) &:& e^{n}_{ij}  &\mapsto &
                    a_{+}(n,i,j) e^{n+\half}_{i-\half ,j-\half }
       + a_{-}(n,i,j)  e^{(n-\half )}_{i-\half ,j-\half },\label{alpha} \\
\lambda_{q}(\beta_{q}) &:& e^{n}_{ij}  &\mapsto & b_{+}(n,i,j)
                 e^{(n+\half )}_{i+\half ,j-\half }
       + b_{-}(n,i,j)  e^{(n-\half )}_{i+\half ,j-\half },\label{beta}
\end{IEEEeqnarray}
%%-------------------------------------
where
%%-------------------------------------
\begin{IEEEeqnarray}{rCl}
a_{+}(n,i,j)  & = &
q^{2n+i+j+1}\frac{g(n-j+1)g(n-i+1)}{g(2n+1)g(2n+2)},
                          \label{aplus}\\
a_{-}(n,i,j)&=&\frac{g(n+j)g(n+i)}{g(2n)g(2n+1)},\label{aminus}\\
b_{+}(n,i,j)&=&
 -q^{n+j}\frac{g(n-j+1)g(n+i+1)}{g(2n+1)g(2n+2)}, \label{bplus}\\
b_{-}(n,i,j)&=&
    q^{n+i}\frac{g(n+j)g(n-i)}{g(2n)g(2n+1)}.\label{bminus}
\end{IEEEeqnarray}
%%-------------------------------------
%%
One can view the representations $\lambda_{q}$ as acting on the single Hilbert
space $\ell^{2}(\Gamma)$ with an orthonormal basis $\{e^{n}_{ij}:
n\in\half\mathbb{N}, i,j=-n,-n+1,\ldots,n-1,n\}$. Thus for each $q\neq 0$, we have
a faithful representation $\pi_{q}$ of $A_{q}$ acting on $\mathcal{H}_{\pi}$ and
another faithful representation $\lambda_{q}$ acting on $\ell^{2}(\Gamma)$. The
actions of $\alpha_{q}$ and $\beta_{q}$ are shown in the following diagram, where
the black double headed arrows represent the action of $\alpha_{q}$ and the red
arrows represent the action of $\beta_{q}$. For both, the solid colored arrow
stands for the second terms in (\ref{alpha}) and (\ref{beta}) and the dashed arrow
stands for the first terms in (\ref{alpha}) and (\ref{beta}).
%%--------------------------------------------------------------------
% \input suq2diagrams-nonzero-q.tex
%%--------------------------------------------------------------------
%%--------------------------------------------------------------------
\begin{center}
\tdplotsetmaincoords{70}{120}
\begin{tikzpicture}[tdplot_main_coords]
	\draw [fill=red] (0,0,3) circle (1.5pt) node[above]{\textcolor{red}{\scriptsize $e^{(0)}_{0,0}$}};
	\draw [fill=blue] (-1,-1,2) circle (1.5pt) node[anchor=south east]{\textcolor{blue}{\scriptsize $e^{(\frac{1}{2})}_{-\frac{1}{2},\frac{1}{2}}$}};
	\draw [fill=blue] (-1,1,2) circle (1.5pt) node[anchor=south west]{\textcolor{blue}{\scriptsize $e^{(\frac{1}{2})}_{\frac{1}{2},\frac{1}{2}}$}};
	\draw [fill=blue] (1,-1,2) circle (1.5pt) node[anchor=south east]{\textcolor{blue}{\scriptsize $e^{(\frac{1}{2})}_{-\frac{1}{2},-\frac{1}{2}}$}};
	\draw [fill=blue] (1,1,2) circle (1.5pt);
	\draw [fill=orange] (-2,-2,1) circle (1.5pt);
	\draw [fill=orange] (-2,0,1) circle (1.5pt);
	\draw [fill=orange] (-2,2,1) circle (1.5pt) node[anchor=south west]{\textcolor{orange}{\scriptsize $e^{(1)}_{1,1}$}};
	\draw [fill=orange] (0,-2,1) circle (1.5pt);
	\draw [fill=orange] (0,0,1) circle (1.5pt);
	\draw [fill=orange] (0,2,1) circle (1.5pt);
	\draw [fill=orange] (2,-2,1) circle (1.5pt) node[anchor=south east]{\textcolor{orange}{\scriptsize $e^{(1)}_{-1,-1}$}};
	\draw [fill=orange] (2,0,1) circle (1.5pt);
	\draw [fill=orange] (2,2,1) circle (1.5pt) node[anchor=north west]{\textcolor{orange}{\scriptsize $e^{(1)}_{1,-1}$}};
	\draw [fill=gray] (-3,-3,0) circle (1.5pt);
	\draw [fill=gray] (-3,-1,0) circle (1.5pt);
	\draw [fill=gray] (-3,1,0) circle (1.5pt);
	\draw [fill=gray] (-3,3,0) circle (1.5pt) node[anchor=south west]{\textcolor{gray}{\scriptsize $e^{(\frac{3}{2})}_{\frac{3}{2},\frac{3}{2}}$}};
	\draw [fill=gray] (-1,-3,0) circle (1.5pt);
	\draw [fill=gray] (-1,-1,0) circle (1.5pt);
	\draw [fill=gray] (-1,1,0) circle (1.5pt);
	\draw [fill=gray] (-1,3,0) circle (1.5pt);
	\draw [fill=gray] (1,-3,0) circle (1.5pt);
	\draw [fill=gray] (1,-1,0) circle (1.5pt);
	\draw [fill=gray] (1,1,0) circle (1.5pt);
	\draw [fill=gray] (1,3,0) circle (1.5pt);
	\draw [fill=gray] (3,-3,0) circle (1.5pt) node[anchor=north east]{\textcolor{gray}{\scriptsize $e^{(\frac{3}{2})}_{-\frac{3}{2},-\frac{3}{2}}$}};
	\draw [fill=gray] (3,-1,0) circle (1.5pt);
	\draw [fill=gray] (3,1,0) circle (1.5pt);
	\draw [fill=gray] (3,3,0) circle (1.5pt) node[anchor=north west]{\textcolor{gray}{\scriptsize $e^{(\frac{3}{2})}_{\frac{3}{2},-\frac{3}{2}}$}};
	\draw[gray] (-3,-3,0) +(0,0,0) -- +(2,0,0) -- + (4,0,0) -- +(6,0,0);
	\draw[gray] (-3,-3,0) +(0,2,0) -- +(2,2,0) -- + (4,2,0) -- +(6,2,0);
	\draw[gray] (-3,-3,0) +(0,4,0) -- +(2,4,0) -- + (4,4,0) -- +(6,4,0);
	\draw[gray] (-3,-3,0) +(0,6,0) -- +(2,6,0) -- + (4,6,0) -- +(6,6,0);
	\draw[gray] (-3,-3,0) +(6,0,0) -- +(6,2,0) -- +(6,4,0) -- +(6,6,0);
	\draw[gray] (-3,-3,0) +(4,0,0) -- +(4,2,0) -- +(4,4,0) -- +(4,6,0);
	\draw[gray] (-3,-3,0) +(2,0,0) -- +(2,2,0) -- +(2,4,0) -- +(2,6,0);
	\draw[gray] (-3,-3,0) +(0,0,0) -- +(0,2,0) -- +(0,4,0) -- +(0,6,0);
	\draw[orange] (-2,-2,1) +(0,0,0) -- +(2,0,0) -- + (4,0,0);
	\draw[orange] (-2,-2,1) +(0,2,0) -- +(2,2,0) -- + (4,2,0);
	\draw[orange] (-2,-2,1) +(0,4,0) -- +(2,4,0) -- + (4,4,0);
	\draw[orange] (-2,-2,1) +(4,0,0) -- +(4,2,0) -- +(4,4,0);
	\draw[orange] (-2,-2,1) +(2,0,0) -- +(2,2,0) -- +(2,4,0);
	\draw[orange] (-2,-2,1) +(0,0,0) -- +(0,2,0) -- +(0,4,0);
	\draw[blue] (-1,-1,2) +(0,0,0) -- +(2,0,0) -- +(2,2,0) -- +(0,2,0) -- +(0,0,0);
	\draw[-{Stealth}{Stealth},black] (0,0,1) to node{} (1,-1,2);
	\draw[-{Stealth}{Stealth},black, dashed] (0,0,1) to node{} (1,-1,0);
	\draw[-{Stealth},red] (0,0,1) to node{} (1,1,2);
	\draw[-{Stealth},red, dashed] (0,0,1) to node{} (1,1,0);
	\draw[] (-1,-1,-2.5)  node[below]{\scriptsize $\lambda_{q} \text{ on }\ell^{2}(\Gamma)$};
\end{tikzpicture}
\end{center}
%%--------------------------------------------------------------------

%%--------------------------------------------------------------------
\subsection{The crystallised $C^{*}$-algebra $A_{0}$}
%%--------------------------------------------------------------------
The relations (\ref{eq:suq2-relations}) at
$q=0$, namely,
%%-------------------------------------
\begin{IEEEeqnarray*}{rClrCl}
\alpha_{0}^*\alpha_{0}+\beta_{0}^*\beta_{0} &=& I,\qquad
       & \alpha_{0}\beta_{0} &=& 0,\\
\alpha_{0}\alpha_{0}^* &=& I,
    &\alpha_{0}\beta_{0}^* &=& 0,\yesnumber\label{eq:su02-relations}\\
&&&\beta_{0}^*\beta_{0} &=& \beta_{0}\beta_{0}^*,
\end{IEEEeqnarray*}
%%-------------------------------------
generate a universal $C^*$-algebra $A_{0}$ which is isomorphic to the
$C^{*}$-algebra $A_{q}\equiv C(SU_{q}(2))$ for $q\neq 0$. Analogous to the
representations $\pi_{q}$ and $\lambda_{q}$ for $A_{q}$, the $C^{*}$-algebra
$A_{0}$ has the following faithful representations on $\mathcal{H}_{\pi}$ and $\mathcal{H}$
respectively:
%%-------------------------------------
\begin{IEEEeqnarray}{rClrCl}
	\pi_{0}(\alpha_{0}) &=& S\otimes I,\qquad
	   & \pi_{0}(\beta_{0}) &=& P_{0}\otimes S.\label{eqn:qzero-1}
\end{IEEEeqnarray}
\begin{IEEEeqnarray}{rClrCl}
	\lambda_{0}(\alpha_{0})e_{ij}^{n} &=&
	  \begin{cases}
	 e_{i-\half,j-\half}^{n-\half} & \text{if }i,j>-n,\qquad\\
	 0 & \text{otherwise.}
	  \end{cases}
	& \lambda_{0}(\beta_{0})e_{ij}^{n} &= &
	  \begin{cases}
	      e_{i-\half,j+\half}^{n+\half} & \text{if $i=-n$},\\
		  -e_{i-\half,j+\half}^{n-\half} & \text{if $j=-n$},\\
		  0 & \text{otherwise.}
	  \end{cases}\label{eqn:qzero-2}
\end{IEEEeqnarray}
%%-------------------------------------

\brmrk
As mentioned in the introduction, Pal \& Giri \cite{GirPal-2022tv} have introduced the notion of crystallisation of quantised function algebras for $q$-deformations of classical compact Lie groups in the type $A_{n}$ case and Matassa \& Yuncken \cite{MatYun-2022aa} for the general case. The $C^*$-algebra described above is the crystallised $C^*$-algebra according to their notions for the $SU_{q}(2)$ case. However, in this case the $C^*$-algebra is very easy to obtain and was described by Woronowicz in \cite{Wor-1987aa}.
\ermrk

%%--------------------------------------------------------------------
\subsection{The unitary}
%%--------------------------------------------------------------------
In this section, we obtain the unitary that will be used in the main theorem for
the equivalence. This is done by passing to $q=0$ and studying the behaviour of
the operators $\lambda_{0}(\alpha_{0})$ and $\lambda_{0}(\beta_{0})$. As observed
by Connes in~\cite{Con-2004ab}, significant simplifications happen at $q=0$: one
term from (\ref{alpha}) and (\ref{beta}) disappear and the actions of $\alpha_{0}$
and $\beta_{0}$ becomes simpler, as shown in the following diagram:
%%--------------------------------------------------------------------
% \input suq2diagrams-zero-q.tex
%%--------------------------------------------------------------------
%%--------------------------------------------------------------------
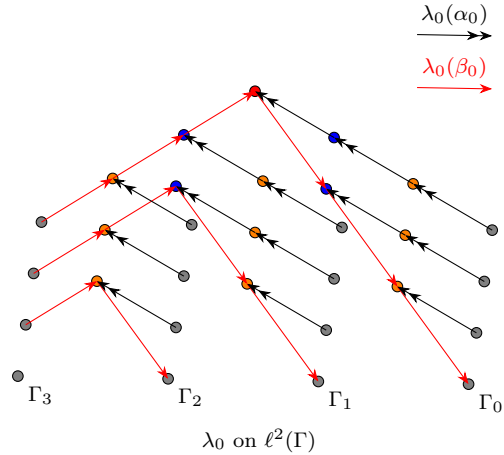
\begin{center}
\tdplotsetmaincoords{70}{93}
\begin{tikzpicture}[tdplot_main_coords]
	\draw [fill=red] (0,0,3) circle (2pt);
	\draw [fill=blue] (-1,-1,2) circle (2pt);
	\draw [fill=blue] (-1,1,2) circle (2pt);
	\draw [fill=blue] (1,-1,2) circle (2pt);
	\draw [fill=blue] (1,1,2) circle (2pt);
	\draw [fill=orange] (-2,-2,1) circle (2pt);
	\draw [fill=orange] (-2,0,1) circle (2pt);
	\draw [fill=orange] (-2,2,1) circle (2pt);
	\draw [fill=orange] (0,-2,1) circle (2pt);
	\draw [fill=orange] (0,0,1) circle (2pt);
	\draw [fill=orange] (0,2,1) circle (2pt);
	\draw [fill=orange] (2,-2,1) circle (2pt);
	\draw [fill=orange] (2,0,1) circle (2pt);
	\draw [fill=orange] (2,2,1) circle (2pt);
	\draw [fill=gray] (-3,-3,0) circle (2pt);
	\draw [fill=gray] (-3,-1,0) circle (2pt);
	\draw [fill=gray] (-3,1,0) circle (2pt);
	\draw [fill=gray] (-3,3,0) circle (2pt);
	\draw [fill=gray] (-1,-3,0) circle (2pt);
	\draw [fill=gray] (-1,-1,0) circle (2pt);
	\draw [fill=gray] (-1,1,0) circle (2pt);
	\draw [fill=gray] (-1,3,0) circle (2pt);
	\draw [fill=gray] (1,-3,0) circle (2pt);
	\draw [fill=gray] (1,-1,0) circle (2pt);
	\draw [fill=gray] (1,1,0) circle (2pt);
	\draw [fill=gray] (1,3,0) circle (2pt);
	\draw [fill=gray] (3,-3,0) circle (2pt) node[anchor=north west]{\textcolor{black}{\scriptsize $\Gamma_{3}$}};
	\draw [fill=gray] (3,-1,0) circle (2pt) node[anchor=north west]{\textcolor{black}{\scriptsize $\Gamma_{2}$}};
	\draw [fill=gray] (3,1,0) circle (2pt) node[anchor=north west]{\textcolor{black}{\scriptsize $\Gamma_{1}$}};
	\draw [fill=gray] (3,3,0) circle (2pt) node[anchor=north west]{\textcolor{black}{\scriptsize $\Gamma_{0}$}};
	\draw[red,-{Stealth}] (-3,-3,0)  to node{}  (-2,-2,1);
	\draw[red,-{Stealth}] (-2,-2,1)  to node{}  (-1,-1,2);
	\draw[red,-{Stealth}] (-1,-1,2)  to node{}  (0,0,3);
	\draw[red,-{Stealth}] (0,0,3)  to node{}  (1,1,2);
	\draw[red,-{Stealth}] (1,1,2)  to node{}  (2,2,1);
	\draw[red,-{Stealth}] (2,2,1)  to node{}  (3,3,0);
	\draw[red,-{Stealth}] (-1,-3,0)  to node{}  (0,-2,1);
	\draw[red,-{Stealth}] (0,-2,1)  to node{}  (1,-1,2);
	\draw[red,-{Stealth}] (1,-1,2)  to node{}  (2,0,1);
	\draw[red,-{Stealth}] (2,0,1)  to node{}  (3,1,0);
	\draw[red,-{Stealth}] (1,-3,0)  to node{}  (2,-2,1);
	\draw[red,-{Stealth}] (2,-2,1)  to node{}  (3,-1,0);
	\draw[-{Stealth}{Stealth}] (-3,-1,0)  to node{}  (-2,-2,1);
	\draw[-{Stealth}{Stealth}] (-3,1,0)  to node{}  (-2,0,1);
	\draw[-{Stealth}{Stealth}] (-2,0,1)  to node{}  (-1,-1,2);
	\draw[-{Stealth}{Stealth}] (-3,3,0)  to node{}  (-2,2,1);
	\draw[-{Stealth}{Stealth}] (-2,2,1)  to node{}  (-1,1,2);
	\draw[-{Stealth}{Stealth}] (-1,1,2)  to node{}  (0,0,3);
	\draw[-{Stealth}{Stealth}] (-1,3,0)  to node{}  (0,2,1);
	\draw[-{Stealth}{Stealth}] (0,2,1)  to node{}  (1,1,2);
	\draw[-{Stealth}{Stealth}] (1,3,0)  to node{}  (2,2,1);
	\draw[-{Stealth}{Stealth}] (-1,-1,0)  to node{}  (0,-2,1);
	\draw[-{Stealth}{Stealth}] (-1,1,0)  to node{}  (0,0,1);
	\draw[-{Stealth}{Stealth}] (0,0,1)  to node{}  (1,-1,2);
	\draw[-{Stealth}{Stealth}] (1,1,0)  to node{}  (2,0,1);
	\draw[-{Stealth}{Stealth}] (1,-1,0)  to node{}  (2,-2,1);
	\draw[-{Stealth}{Stealth}] (-3,2,2.75)  to node[above]{\scriptsize $\lambda_{0}(\alpha_{0})$}  (-3,3,2.75);
	\draw[red,-{Stealth}] (-3,2,2)  to node[above]{\scriptsize $\lambda_{0}(\beta_{0})$}  (-3,3,2);
	\draw[] (-1,0,-2)  node[below]{\scriptsize $\lambda_{0} \text{ on }\ell^{2}(\Gamma)$};
\end{tikzpicture}
\end{center}
%%--------------------------------------------------------------------
%%--------------------------------------------------------------------
Let us denote the sheet consisting of the right and rear face of the pyramid by
$\Gamma_{0}$, i.e.
%%-------------------------------------
\[
\Gamma_{0}=\left\{(n,i,j): n\in  \half\mathbb{N},
    \; i,j\in \{-n, -n+1,\ldots,n-1,n\},\; \max\{i,j\}=n\right\}.
\]
%%-------------------------------------
Then one can naturally identify $\Gamma_{0}$ with $\mathbb{N}\times \mathbb{Z}$ as
the next diagram illustrates, which means there is a unitary $V_{0}$ between
$\ell^{2}(\Gamma_{0})$ and $\ell^{2}(\mathbb{N}\times\mathbb{Z})$. As the actions
of $\lambda_{0}(\alpha_{0})$ and $\lambda_{0}(\beta_{0})$ keep the sheet
$\Gamma_{0}$ invariant, through this unitary, they are equivalent to operators on
$\ell^{2}(\mathbb{N}\times\mathbb{Z})$. In fact they turn out to be precisely
$\pi_{0}(\alpha_{0})$ and $\pi_{0}(\beta_{0})$. Next, taking the sheet $\Gamma_{0}$
out, $\Gamma$ is left with an identical replica of itself. Denote by $\Gamma_{1}$
by the union of the right and rear face of this remaining part. Continue in this
fashion and define, for $k\in\mathbb{N}$,
%%-------------------------------------
\[
\Gamma_{k}=\left\{(n,i,j): n\in \half\mathbb{N}, \; n\geq \frac{k}{2},
   \; i,j\in \{-n, -n+1,\ldots,n-1,n\},\; \max\{i,j\}=n\right\}.
\]
%%-------------------------------------
Then just like $\Gamma_{0}$, each $\Gamma_{k}$ is kept invariant by
$\lambda_{0}(\alpha_{0})$ and $\lambda_{0}(\beta_{0})$, can naturally be identified
with $\mathbb{N}\times \mathbb{Z}$ through a unitary $V_{k}$ and
$V_{k}\lambda_{0}(\alpha_{0})V_{k}^{*}$ and $V_{k}\lambda_{0}(\beta_{0})V_{k}^{*}$
are precisely the operators $\pi_{0}(\alpha_{0})$ and $\pi_{0}(\beta_{0})$ on
$\ell^{2}(\mathbb{N}\times\mathbb{Z})$.\\
%%--------------------------------------------------------------------
% \input suq2diagrams-unitary-1.tex
%%--------------------------------------------------------------------
%%--------------------------------------------------------------------
\tdplotsetmaincoords{70}{93}
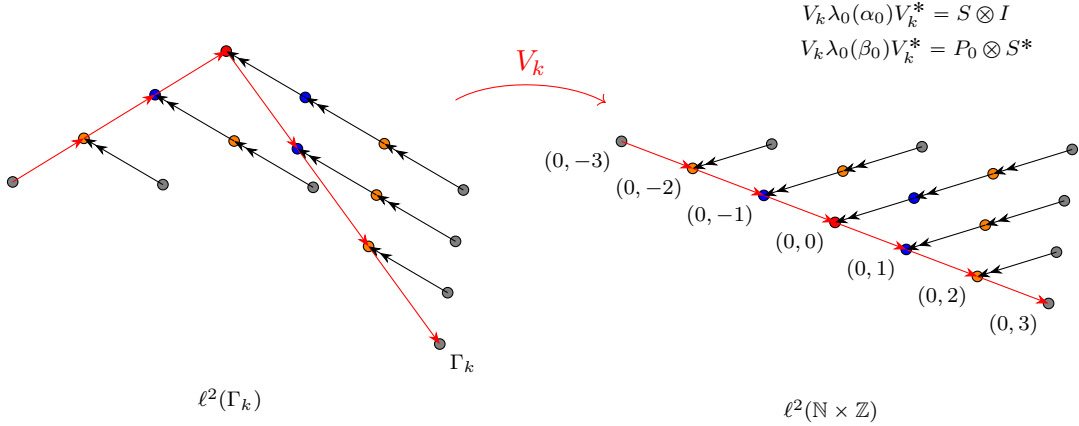
\begin{tikzpicture}[tdplot_main_coords]
	\draw [fill=red] (0,0,3) circle (2pt);
	\draw [fill=blue] (-1,-1,2) circle (2pt);
	\draw [fill=blue] (-1,1,2) circle (2pt);
	\draw [fill=blue] (1,1,2) circle (2pt);
	\draw [fill=orange] (-2,-2,1) circle (2pt);
	\draw [fill=orange] (-2,0,1) circle (2pt);
	\draw [fill=orange] (-2,2,1) circle (2pt);
	\draw [fill=orange] (0,2,1) circle (2pt);
	\draw [fill=orange] (2,2,1) circle (2pt);
	\draw [fill=gray] (-3,-3,0) circle (2pt);
	\draw [fill=gray] (-3,-1,0) circle (2pt);
	\draw [fill=gray] (-3,1,0) circle (2pt);
	\draw [fill=gray] (-3,3,0) circle (2pt);
	\draw [fill=gray] (-1,3,0) circle (2pt);
	\draw [fill=gray] (1,3,0) circle (2pt);
	\draw [fill=gray] (3,3,0) circle (2pt) node[anchor=north west]{\textcolor{black}{\scriptsize $\Gamma_{k}$}};
	\draw[red,-{Stealth}] (-3,-3,0) -- (-2,-2,1) node[anchor=north east]{};
	\draw[red,-{Stealth}] (-2,-2,1) -- (-1,-1,2) node[anchor=north east]{};
	\draw[red,-{Stealth}] (-1,-1,2) -- (0,0,3) node[anchor=north east]{};
	\draw[red,-{Stealth}] (0,0,3) -- (1,1,2) node[anchor=north east]{};
	\draw[red,-{Stealth}] (1,1,2) -- (2,2,1) node[anchor=north east]{};
	\draw[red,-{Stealth}] (2,2,1) -- (3,3,0) node[anchor=north east]{};
	\draw[-{Stealth}{Stealth}] (-3,-1,0) -- (-2,-2,1) node[anchor=north east]{};
	\draw[-{Stealth}{Stealth}] (-3,1,0) -- (-2,0,1) node[anchor=north east]{};
	\draw[-{Stealth}{Stealth}] (-2,0,1) -- (-1,-1,2) node[anchor=north east]{};
	\draw[-{Stealth}{Stealth}] (-3,3,0) -- (-2,2,1) node[anchor=north east]{};
	\draw[-{Stealth}{Stealth}] (-2,2,1) -- (-1,1,2) node[anchor=north east]{};
	\draw[-{Stealth}{Stealth}] (-1,1,2) -- (0,0,3) node[anchor=north east]{};
	\draw[-{Stealth}{Stealth}] (-1,3,0) -- (0,2,1) node[anchor=north east]{};
	\draw[-{Stealth}{Stealth}] (0,2,1) -- (1,1,2) node[anchor=north east]{};
	\draw[-{Stealth}{Stealth}] (1,3,0) -- (2,2,1) node[anchor=north east]{};
	\draw [->,out=30,in=150,looseness=0.75,red] (-1,3,2) to node[above]{$V_{k}$}  (-1,5,2);
	\draw [fill=red] (-2,8,0) circle (2pt) node[anchor=north east]{\scriptsize $(0,0)$};
	\draw [fill=blue] (-3,7,0) circle (2pt) node[anchor=north east]{\scriptsize $(0,-1)$};
	\draw [fill=blue] (-3,9,0) circle (2pt);
	\draw [fill=blue] (-1,9,0) circle (2pt) node[anchor=north east]{\scriptsize $(0,1)$};
	\draw [fill=orange] (-4,6,0) circle (2pt) node[anchor=north east]{\scriptsize $(0,-2)$};
	\draw [fill=orange] (-4,8,0) circle (2pt);
	\draw [fill=orange] (-4,10,0) circle (2pt);
	\draw [fill=orange] (-2,10,0) circle (2pt);
	\draw [fill=orange] (0,10,0) circle (2pt) node[anchor=north east]{\scriptsize $(0,2)$};
	\draw [fill=gray] (-5,5,0) circle (2pt) node[anchor=north east]{\scriptsize $(0,-3)$};
	\draw [fill=gray] (-5,7,0) circle (2pt);
	\draw [fill=gray] (-5,9,0) circle (2pt);
	\draw [fill=gray] (-5,11,0) circle (2pt);
	\draw [fill=gray] (-3,11,0) circle (2pt);
	\draw [fill=gray] (-1,11,0) circle (2pt);
	\draw [fill=gray] (1,11,0) circle (2pt) node[anchor=north east]{\scriptsize $(0,3)$};
	\draw[red,-{Stealth}] (-5,5,0) -- (-4,6,0) node[anchor=north east]{};
	\draw[red,-{Stealth}] (-4,6,0) -- (-3,7,0) node[anchor=north east]{};
	\draw[red,-{Stealth}] (-3,7,0) -- (-2,8,0) node[anchor=north east]{};
	\draw[red,-{Stealth}] (-2,8,0) -- (-1,9,0) node[anchor=north east]{};
	\draw[red,-{Stealth}] (-1,9,0) -- (0,10,0) node[anchor=north east]{};
	\draw[red,-{Stealth}] (0,10,0) -- (1,11,0) node[anchor=north east]{};
	\draw[-{Stealth}{Stealth}] (-5,7,0) -- (-4,6,0) node[anchor=north east]{};
	\draw[-{Stealth}{Stealth}] (-5,9,0) -- (-4,8,0) node[anchor=north east]{};
	\draw[-{Stealth}{Stealth}] (-4,8,0) -- (-3,7,0) node[anchor=north east]{};
	\draw[-{Stealth}{Stealth}] (-5,11,0) -- (-4,10,0) node[anchor=north east]{};
	\draw[-{Stealth}{Stealth}] (-4,10,0) -- (-3,9,0) node[anchor=north east]{};
	\draw[-{Stealth}{Stealth}] (-3,9,0) -- (-2,8,0) node[anchor=south west]{};
	\draw[-{Stealth}{Stealth}] (-3,11,0) -- (-2,10,0) node[anchor=north east]{};
	\draw[-{Stealth}{Stealth}] (-2,10,0) -- (-1,9,0) node[anchor=north east]{};
	\draw[-{Stealth}{Stealth}] (-1,11,0) -- (0,10,0) node[anchor=north east]{};
	\draw[] (-1,9,3)  node[above]{\scriptsize $V_{k}\lambda_{0}(\alpha_{0})V_{k}^{*}=S\otimes I$};
	\draw[] (-1,9,2.5)  node[above]{\scriptsize \quad$V_{k}\lambda_{0}(\beta_{0})V_{k}^{*}=P_{0}\otimes S^{*}$};
	\draw[] (-1,0,-2)  node[below]{\scriptsize $\ell^{2}(\Gamma_{k})$};
	\draw[] (-1,8,-2)  node[below]{\scriptsize $\ell^{2}(\mathbb{N}\times\mathbb{Z})$};
\end{tikzpicture}
%%--------------------------------------------------------------------

\noindent
Taking the direct sum of all these unitaries $V_{k}$, one gets a unitary $U$ from
$\ell^{2}(\Gamma)$ to the space $\ell^{2}(\mathbb{N})\otimes
\ell^{2}(\mathbb{N}\times\mathbb{Z})$:\\
%%--------------------------------------------------------------------
% \input suq2diagrams-unitary-2.tex
%%--------------------------------------------------------------------
%%--------------------------------------------------------------------
\tdplotsetmaincoords{70}{93}
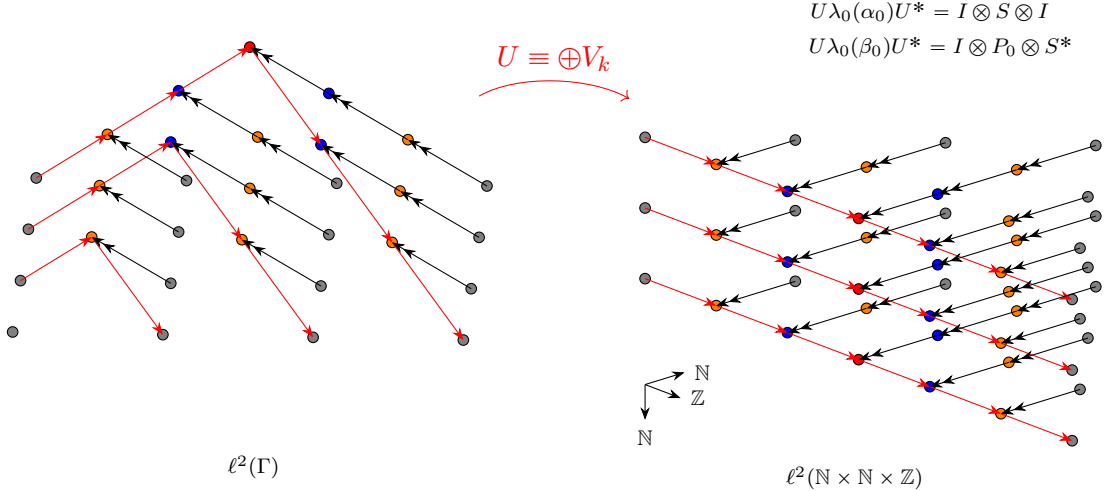
\begin{tikzpicture}[tdplot_main_coords]
	\draw [fill=red] (0,0,3) circle (2pt);
	\draw [fill=blue] (-1,-1,2) circle (2pt);
	\draw [fill=blue] (-1,1,2) circle (2pt);
	\draw [fill=blue] (1,-1,2) circle (2pt);
	\draw [fill=blue] (1,1,2) circle (2pt);
	\draw [fill=orange] (-2,-2,1) circle (2pt);
	\draw [fill=orange] (-2,0,1) circle (2pt);
	\draw [fill=orange] (-2,2,1) circle (2pt);
	\draw [fill=orange] (0,-2,1) circle (2pt);
	\draw [fill=orange] (0,0,1) circle (2pt);
	\draw [fill=orange] (0,2,1) circle (2pt);
	\draw [fill=orange] (2,-2,1) circle (2pt);
	\draw [fill=orange] (2,0,1) circle (2pt);
	\draw [fill=orange] (2,2,1) circle (2pt);
	\draw [fill=gray] (-3,-3,0) circle (2pt);
	\draw [fill=gray] (-3,-1,0) circle (2pt);
	\draw [fill=gray] (-3,1,0) circle (2pt);
	\draw [fill=gray] (-3,3,0) circle (2pt);
	\draw [fill=gray] (-1,-3,0) circle (2pt);
	\draw [fill=gray] (-1,-1,0) circle (2pt);
	\draw [fill=gray] (-1,1,0) circle (2pt);
	\draw [fill=gray] (-1,3,0) circle (2pt);
	\draw [fill=gray] (1,-3,0) circle (2pt);
	\draw [fill=gray] (1,-1,0) circle (2pt);
	\draw [fill=gray] (1,1,0) circle (2pt);
	\draw [fill=gray] (1,3,0) circle (2pt);
	\draw [fill=gray] (3,-3,0) circle (2pt);
	\draw [fill=gray] (3,-1,0) circle (2pt);
	\draw [fill=gray] (3,1,0) circle (2pt);
	\draw [fill=gray] (3,3,0) circle (2pt);
	\draw[red,-{Stealth}] (-3,-3,0)   to node{}   (-2,-2,1);
	\draw[red,-{Stealth}] (-2,-2,1)   to node{}   (-1,-1,2);
	\draw[red,-{Stealth}] (-1,-1,2)   to node{}   (0,0,3);
	\draw[red,-{Stealth}] (0,0,3)   to node{}   (1,1,2);
	\draw[red,-{Stealth}] (1,1,2)   to node{}   (2,2,1);
	\draw[red,-{Stealth}] (2,2,1)   to node{}   (3,3,0);
	\draw[red,-{Stealth}] (-1,-3,0)   to node{}   (0,-2,1);
	\draw[red,-{Stealth}] (0,-2,1)   to node{}   (1,-1,2);
	\draw[red,-{Stealth}] (1,-1,2)   to node{}   (2,0,1);
	\draw[red,-{Stealth}] (2,0,1)   to node{}   (3,1,0);
	\draw[red,-{Stealth}] (1,-3,0)   to node{}   (2,-2,1);
	\draw[red,-{Stealth}] (2,-2,1)   to node{}   (3,-1,0);
	\draw[-{Stealth}{Stealth}] (-3,-1,0)   to node{}   (-2,-2,1);
	\draw[-{Stealth}{Stealth}] (-3,1,0)   to node{}   (-2,0,1);
	\draw[-{Stealth}{Stealth}] (-2,0,1)   to node{}   (-1,-1,2);
	\draw[-{Stealth}{Stealth}] (-3,3,0)   to node{}   (-2,2,1);
	\draw[-{Stealth}{Stealth}] (-2,2,1)   to node{}   (-1,1,2);
	\draw[-{Stealth}{Stealth}] (-1,1,2)   to node{}   (0,0,3);
	\draw[-{Stealth}{Stealth}] (-1,3,0)   to node{}   (0,2,1);
	\draw[-{Stealth}{Stealth}] (0,2,1)   to node{}   (1,1,2);
	\draw[-{Stealth}{Stealth}] (1,3,0)   to node{}   (2,2,1);
	\draw[-{Stealth}{Stealth}] (-1,-1,0)   to node{}   (0,-2,1);
	\draw[-{Stealth}{Stealth}] (-1,1,0)   to node{}   (0,0,1);
	\draw[-{Stealth}{Stealth}] (0,0,1)  to node{}   (1,-1,2);
	\draw[-{Stealth}{Stealth}] (1,1,0)  to node{}  (2,0,1);
	\draw[-{Stealth}{Stealth}] (1,-1,0)  to node{}  (2,-2,1);
	\draw [->,out=30,in=150,looseness=0.75,red] (-1,3,2) to node[above]{$U\equiv\oplus V_{k}$}  (-1,5,2);
	\draw [fill=red] (-2,8,0) circle (2pt);
	\draw [fill=blue] (-3,7,0) circle (2pt);
	\draw [fill=blue] (-3,9,0) circle (2pt);
	\draw [fill=blue] (-1,9,0) circle (2pt) ;
	\draw [fill=orange] (-4,6,0) circle (2pt);
	\draw [fill=orange] (-4,8,0) circle (2pt);
	\draw [fill=orange] (-4,10,0) circle (2pt);
	\draw [fill=orange] (-2,10,0) circle (2pt);
	\draw [fill=orange] (0,10,0) circle (2pt) ;
	\draw [fill=gray] (-5,5,0) circle (2pt) ;
	\draw [fill=gray] (-5,7,0) circle (2pt);
	\draw [fill=gray] (-5,9,0) circle (2pt);
	\draw [fill=gray] (-5,11,0) circle (2pt);
	\draw [fill=gray] (-3,11,0) circle (2pt);
	\draw [fill=gray] (-1,11,0) circle (2pt);
	\draw [fill=gray] (1,11,0) circle (2pt) ;
	\draw[red,-{Stealth}] (-5,5,0) -- (-4,6,0);
	\draw[red,-{Stealth}] (-4,6,0) -- (-3,7,0) ;
	\draw[red,-{Stealth}] (-3,7,0) -- (-2,8,0) ;
	\draw[red,-{Stealth}] (-2,8,0) -- (-1,9,0) ;
	\draw[red,-{Stealth}] (-1,9,0) -- (0,10,0) ;
	\draw[red,-{Stealth}] (0,10,0) -- (1,11,0) ;
	\draw[-{Stealth}{Stealth}] (-5,7,0) -- (-4,6,0) ;
	\draw[-{Stealth}{Stealth}] (-5,9,0) -- (-4,8,0) ;
	\draw[-{Stealth}{Stealth}] (-4,8,0) -- (-3,7,0) ;
	\draw[-{Stealth}{Stealth}] (-5,11,0) -- (-4,10,0) ;
	\draw[-{Stealth}{Stealth}] (-4,10,0) -- (-3,9,0) ;
	\draw[-{Stealth}{Stealth}] (-3,9,0) -- (-2,8,0) ;
	\draw[-{Stealth}{Stealth}] (-3,11,0) -- (-2,10,0) ;
	\draw[-{Stealth}{Stealth}] (-2,10,0) -- (-1,9,0) ;
	\draw[-{Stealth}{Stealth}] (-1,11,0) -- (0,10,0) ;
	\draw[] (-1,9,3)  node[above]{\scriptsize $U\lambda_{0}(\alpha_{0})U^{*}=I\otimes S\otimes I$};
	\draw[] (-1,9,2.5)  node[above]{\scriptsize $\quad U\lambda_{0}(\beta_{0})U^{*}=I\otimes P_{0}\otimes S^{*}$};
	\draw[] (-1,0,-3)  node[below]{\scriptsize $\ell^{2}(\Gamma)$};
	\draw[] (-1,8,-3)  node[below]{\scriptsize $\ell^{2}(\mathbb{N}\times\mathbb{N}\times\mathbb{Z})$};
	\draw [fill=red] (-2,8, -1) circle (2pt);
	\draw [fill=blue] (-3,7, -1) circle (2pt);
	\draw [fill=blue] (-3,9, -1) circle (2pt);
	\draw [fill=blue] (-1,9, -1) circle (2pt) ;
	\draw [fill=orange] (-4,6, -1) circle (2pt);
	\draw [fill=orange] (-4,8, -1) circle (2pt);
	\draw [fill=orange] (-4,10, -1) circle (2pt);
	\draw [fill=orange] (-2,10, -1) circle (2pt);
	\draw [fill=orange] (0,10, -1) circle (2pt) ;
	\draw [fill=gray] (-5,5, -1) circle (2pt) ;
	\draw [fill=gray] (-5,7, -1) circle (2pt);
	\draw [fill=gray] (-5,9, -1) circle (2pt);
	\draw [fill=gray] (-5,11, -1) circle (2pt);
	\draw [fill=gray] (-3,11, -1) circle (2pt);
	\draw [fill=gray] (-1,11, -1) circle (2pt);
	\draw [fill=gray] (1,11, -1) circle (2pt) ;
	\draw[red,-{Stealth}] (-5,5, -1) -- (-4,6, -1);
	\draw[red,-{Stealth}] (-4,6, -1) -- (-3,7, -1) ;
	\draw[red,-{Stealth}] (-3,7, -1) -- (-2,8, -1) ;
	\draw[red,-{Stealth}] (-2,8, -1) -- (-1,9, -1) ;
	\draw[red,-{Stealth}] (-1,9, -1) -- (0,10, -1) ;
	\draw[red,-{Stealth}] (0,10, -1) -- (1,11, -1) ;
	\draw[-{Stealth}{Stealth}] (-5,7, -1) -- (-4,6, -1) ;
	\draw[-{Stealth}{Stealth}] (-5,9, -1) -- (-4,8, -1) ;
	\draw[-{Stealth}{Stealth}] (-4,8, -1) -- (-3,7, -1) ;
	\draw[-{Stealth}{Stealth}] (-5,11, -1) -- (-4,10, -1) ;
	\draw[-{Stealth}{Stealth}] (-4,10, -1) -- (-3,9, -1) ;
	\draw[-{Stealth}{Stealth}] (-3,9, -1) -- (-2,8, -1) ;
	\draw[-{Stealth}{Stealth}] (-3,11, -1) -- (-2,10, -1) ;
	\draw[-{Stealth}{Stealth}] (-2,10, -1) -- (-1,9, -1) ;
	\draw[-{Stealth}{Stealth}] (-1,11, -1) -- (0,10, -1) ;
	\draw [fill=red] (-2,8,  -2) circle (2pt);
	\draw [fill=blue] (-3,7,  -2) circle (2pt);
	\draw [fill=blue] (-3,9,  -2) circle (2pt);
	\draw [fill=blue] (-1,9,  -2) circle (2pt) ;
	\draw [fill=orange] (-4,6,  -2) circle (2pt);
	\draw [fill=orange] (-4,8,  -2) circle (2pt);
	\draw [fill=orange] (-4,10,  -2) circle (2pt);
	\draw [fill=orange] (-2,10,  -2) circle (2pt);
	\draw [fill=orange] (0,10,  -2) circle (2pt) ;
	\draw [fill=gray] (-5,5,  -2) circle (2pt) ;
	\draw [fill=gray] (-5,7,  -2) circle (2pt);
	\draw [fill=gray] (-5,9,  -2) circle (2pt);
	\draw [fill=gray] (-5,11,  -2) circle (2pt);
	\draw [fill=gray] (-3,11,  -2) circle (2pt);
	\draw [fill=gray] (-1,11,  -2) circle (2pt);
	\draw [fill=gray] (1,11,  -2) circle (2pt) ;
	\draw[red,-{Stealth}] (-5,5,  -2) -- (-4,6,  -2);
	\draw[red,-{Stealth}] (-4,6,  -2) -- (-3,7,  -2) ;
	\draw[red,-{Stealth}] (-3,7,  -2) -- (-2,8,  -2) ;
	\draw[red,-{Stealth}] (-2,8,  -2) -- (-1,9,  -2) ;
	\draw[red,-{Stealth}] (-1,9,  -2) -- (0,10,  -2) ;
	\draw[red,-{Stealth}] (0,10,  -2) -- (1,11,  -2) ;
	\draw[-{Stealth}{Stealth}] (-5,7,  -2) -- (-4,6,  -2) ;
	\draw[-{Stealth}{Stealth}] (-5,9,  -2) -- (-4,8,  -2) ;
	\draw[-{Stealth}{Stealth}] (-4,8,  -2) -- (-3,7,  -2) ;
	\draw[-{Stealth}{Stealth}] (-5,11,  -2) -- (-4,10,  -2) ;
	\draw[-{Stealth}{Stealth}] (-4,10,  -2) -- (-3,9,  -2) ;
	\draw[-{Stealth}{Stealth}] (-3,9,  -2) -- (-2,8,  -2) ;
	\draw[-{Stealth}{Stealth}] (-3,11,  -2) -- (-2,10,  -2) ;
	\draw[-{Stealth}{Stealth}] (-2,10,  -2) -- (-1,9,  -2) ;
	\draw[-{Stealth}{Stealth}] (-1,11,  -2) -- (0,10, -2) ;
	\draw[-{Stealth}]  (-5,5,  -3.5) -- (-4.5,5.5, -3.5);
	\draw[-{Stealth}] (-5,5,  -3.5)  --  (-5.5,5.5, -3.5) ;
	\draw[-{Stealth}] (-5,5,-3.5)  --  (-5,5, -4) ;
	\draw[] (-4.5,5.5, -3.5) node[right]{\scriptsize $\mathbb{Z}$};
	\draw[] (-5.5,5.5, -3.5) node[right]{\scriptsize $\mathbb{N}$};
	\draw[] (-5,5, -4) node[below]{\scriptsize $\mathbb{N}$};
\end{tikzpicture}
%%--------------------------------------------------------------------

\noindent
In the remaining part of this section, we carry out the proof outlined above.

We will denote by $\mathcal{H}_{mult}$ the multiplicity space
$\ell^{2}(\mathbb{N})$ in what follows. The unitary $U:\ell^{2}(\Gamma)\cong
L_2(SU_q(2))\to \mathcal{H}_{mult}\otimes\mathcal{H}_{\pi}$ is given by
%%-------------------------------------
\begin{IEEEeqnarray}{rCl}
	U e_{ij}^{n} &= &
	\begin{cases}
	   e( n-j, n+i,j-i) & \text{if } i<j,\\
	 (-1)^{i-j}e(n-i, n+j, j-i) & \text{if }i\geq j
	\end{cases}\nonumber\\
		&= & (-1)^{(i\vee j) -j} e\bigl(n-(i\vee j), n+(i\wedge j), j-i\bigr).
		   \label{eqn:qzero-3}
\end{IEEEeqnarray}
%%-------------------------------------

We then have
%%-------------------------------------
\begin{IEEEeqnarray}{rCl}
 U^{*}e(r,s,t) &=&
 \begin{cases}
	  e^{\halfof{r+s+|t|}}_{s-\halfof{r+s+|t|}, \halfof{r+s+|t|}+r}
	           & \text{if } t>0,\\
 (-1)^{|t|}e^{\halfof{r+s+|t|}}_{\halfof{r+s+|t|}+r,s-\halfof{r+s+|t|}}
                & \text{if } t\leq 0
					  \end{cases}\nonumber\\
	&=& (-1)^{t_{-}}
	   e^{\halfof{r+s+|t|}}_{\halfof{-r+s-t},\halfof{-r+s+t}}.
	     \label{eqn:qzero-4}
\end{IEEEeqnarray}
%%-------------------------------------

%%--------------------------------------------------------------------
\bthm
For any $a\in A_{0}$, one has $U\lambda_{0}(a)U^{*}=I\otimes \pi_{0}(a)$.
\ethm
%%--------------------------------------------------------------------
\prf
We will show that the equality holds for the generating elements $\alpha_{0}$ and
$\beta_{0}$. From (\ref{eqn:qzero-4}), one has
%%-------------------------------------
\[
	U\lambda_{0}(\alpha_{0})U^{*}e(r,s,t) =
	  (-1)^{t_{-}}U\lambda_{0}(\alpha_{0})
	   e^{\halfof{r}+\halfof{s}+\halfof{|t|}}_{
	      -\halfof{r}+\halfof{s}-\halfof{t},
		    -\halfof{r}+\halfof{s}+\halfof{t}}.
\]
%%-------------------------------------
Note that
\begin{enumerate}
	\item
	 $-\halfof{r}-\halfof{s}-\halfof{|t|}=-\halfof{r}+\halfof{s}-\halfof{t}$
	 if and only if $s+t_{-}= 0$,
	 \item
	 $-\halfof{r}-\halfof{s}-\halfof{|t|}=-\halfof{r}+\halfof{s}+\halfof{t}$
	 if and only if $s+t_{+}=0$,
	 \item
	 $(-\halfof{r}+\halfof{s}-\halfof{t}-\half)
	  \vee(-\halfof{r}+\halfof{s}+\halfof{t}-\half)=
	    -\halfof{r}+\halfof{s}+\halfof{|t|}-\half$,
	 \item
	 $(-\halfof{r}+\halfof{s}-\halfof{t}-\half)
	  \wedge(-\halfof{r}+\halfof{s}+\halfof{t}-\half)=
	  -\halfof{r}+\halfof{s}-\halfof{|t|}-\half$.
\end{enumerate}
Therefore it follows that
%%-------------------------------------
\begin{IEEEeqnarray*}{rCl}
	\IEEEeqnarraymulticol{3}{l}{
	U\lambda_{0}(\alpha_{0})U^{*}e(r,s,t)}\\
	 &=& \begin{cases}
	   (-1)^{t_{-}}Ue^{\halfof{r}+\halfof{s}+\halfof{|t|}-\half}
			   _{-\halfof{r}+\halfof{s}-\halfof{t}-\half,
	 		    -\halfof{r}+\halfof{s}+\halfof{t}-\half}
				 &\text{if } s\neq -t_{-} \text{ and }s\neq -t_{+},\\
	 		   0 &\text{otherwise}
	 		   \end{cases}\\
	&=& \begin{cases} (-1)^{t_{-}+(-\halfof{r}+\halfof{s}+\halfof{|t|}-\half) -
	      (-\halfof{r}+\halfof{s}+\halfof{t}-\half) }&\\
	e\Bigl((\halfof{r}+\halfof{s}+\halfof{|t|}-\half)-
	       (-\halfof{r}+\halfof{s}+\halfof{|t|}-\half),&\\
	\quad (\halfof{r}+\halfof{s}+\halfof{|t|}-\half)+(-\halfof{r}+\halfof{s}-\halfof{|t|}-\half),\;t\Bigr)
	  & \text{if } s\neq -t_{-} \text{ and }s\neq -t_{+},\\
	  0& \text{otherwise}
	 		   \end{cases}\\
	&=&\begin{cases}
	    e(r,s-1,t) & \text{if } t\geq 0 \text{ and }  s\neq 0,\\
		e(r,s-1,t) & \text{if } t< 0  \text{ and } s\neq 0,\\
		0  & \text{otherwise}
	 		   \end{cases}\\
	&=&  (I\otimes \pi_{0}(\alpha_{0}))e(r,s,t).
\end{IEEEeqnarray*}
%%-------------------------------------
Similarly one has
%%-------------------------------------
\[
	U\lambda_{0}(\beta_{0})U^{*}e(r,s,t) =
	  (-1)^{t_{-}}U\lambda_{0}(\beta_{0})
	   e^{\halfof{r}+\halfof{s}+\halfof{|t|}}_{
	      -\halfof{r}+\halfof{s}-\halfof{t},
		    -\halfof{r}+\halfof{s}+\halfof{t}}.
\]
%%-------------------------------------
Note that
\begin{enumerate}
	\item
	 $-\halfof{r}-\halfof{s}-\halfof{|t|}=-\halfof{r}+\halfof{s}-\halfof{t}$
	 if and only if $s+t_{-}= 0$,
	 \item
	 $-\halfof{r}-\halfof{s}-\halfof{|t|}=-\halfof{r}+\halfof{s}+\halfof{t}$
	 if and only if $s+t_{+}=0$.
\end{enumerate}
Therefore
%%-------------------------------------
\begin{IEEEeqnarray*}{rCl}
	U\lambda_{0}(\beta_{0})U^{*}e(r,s,t)
	 &=& \begin{cases}
	   -(-1)^{t_{-}}Ue^{\halfof{r}+\halfof{s}+\halfof{|t|}+\half}
			   _{-\halfof{r}+\halfof{s}-\halfof{t}+\half,
	 		    -\halfof{r}+\halfof{s}+\halfof{t}-\half}
				 &\text{if } s = -t_{+},\\
	   (-1)^{t_{-}}Ue^{\halfof{r}+\halfof{s}+\halfof{|t|}-\half}
			   _{-\halfof{r}+\halfof{s}-\halfof{t}+\half,
	 		    -\halfof{r}+\halfof{s}+\halfof{t}-\half}
				 &\text{if } s \neq -t_{+} \text{ and } s=-t_{-},\\
	 		   0 &\text{ if } s\neq -t_{+} \text{ and } s\neq -t_{-}.
	 		   \end{cases}
\end{IEEEeqnarray*}
%%-------------------------------------
Observe that
\begin{enumerate}
	\item
	one has $s\neq -t_{-}$ and $s\neq -t_{+}$ if and only if
	$s\neq 0$,
	% In this case, $U\lambda_{0}(\beta_{0})U^{*}e^{r}_{s,t}=0$.
	\item
	one has $s \neq -t_{+}$  and  $s=-t_{-}$ if and only if $s=0$ and $t>0$,
	\item
	$s=-t_{+}$ if and only if $s=0$ and $t\leq 0$,
	 \item
	 $(-\halfof{r}+\halfof{s}-\halfof{t}+\half)
	  \vee(-\halfof{r}+\halfof{s}+\halfof{t}-\half)=
	    -\halfof{r}+\halfof{s}+(\halfof{t-1})_{+}$,
	 \item
	 $(-\halfof{r}+\halfof{s}-\halfof{t}+\half)
	  \wedge(-\halfof{r}+\halfof{s}+\halfof{t}-\half)=
	  -\halfof{r}+\halfof{s}-(\halfof{t-1})_{-}$.
\end{enumerate}
Therefore
%%-------------------------------------
\begin{IEEEeqnarray*}{rCl}
	\IEEEeqnarraymulticol{3}{l}{
	U\lambda_{0}(\beta_{0})U^{*}e(r,s,t)}\\
	&=&
	\begin{cases}
		0 & \text{if } s\neq 0,\\
		&\\
		(-1)^{t_{-}+(-\halfof{r}+\halfof{s}+(\halfof{t-1})_{+}) - (-\halfof{r}+\halfof{s}+\halfof{t}-\half) }&\\
			e\Bigl((\halfof{r}+\halfof{s}+\halfof{|t|}-\half)-
			(-\halfof{r}+\halfof{s}+(\halfof{t-1})_{+}),&\\
			\qquad\quad (\halfof{r}+\halfof{s}+\halfof{|t|}-\half)+
			(-\halfof{r}+\halfof{s}-(\halfof{t-1})_{-}),\;t-1\Bigr)
			  & \text{if } s=0 \text{ and }t>0,\\
		&\\
		-(-1)^{t_{-}+(-\halfof{r}+\halfof{s}+(\halfof{t-1})_{+}) - (-\halfof{r}+\halfof{s}+\halfof{t}-\half) }&\\
			e\Bigl((\halfof{r}+\halfof{s}+\halfof{|t|}+\half)-
			(-\halfof{r}+\halfof{s}+(\halfof{t-1})_{+}),&\\
			\qquad\quad \;(\halfof{r}+\halfof{s}+\halfof{|t|}+\half)+
			(-\halfof{r}+\halfof{s}-(\halfof{t-1})_{-}),\;t-1\Bigr)
			  & \text{if } s=0 \text{ and }t\leq 0
	  \end{cases}\\
	&=&
	\begin{cases}
		0 & \text{if } s\neq 0,\\
		e(r,s,t-1) & \text{if } s=0
	\end{cases}\\
	&=& (I\otimes \pi_{0}(\beta_{0}))e(r,s,t).
\end{IEEEeqnarray*}
%%-------------------------------------
Since we have $U\lambda_{0}(a)U^{*}=(I\otimes \pi_{0}(a))$ for $a=\alpha_{0}$ and
$a=\beta_{0}$, the result follows.
\qed

%%--------------------------------------------------------------------
\section{Approximate equivalence for nonzero $q$}
%%--------------------------------------------------------------------
We now come to the main result of the paper, which says that the unitary that gave
us equivalence in the $q=0$ case gives an approximate unitary equivalence between
$\lambda_{q}$ and $I\otimes\pi_{q}$ for $q\neq 0$.
%%--------------------------------------------------------------------
\bthm\label{thm:qnzero-3}
For any $a\in C(SU_q(2))$, one has
\begin{equation}\label{eqn:qnzero-17}
U\lambda_{q}(a)U^{*}-I\otimes \pi_{q}(a)\in \mathscr{T}\otimes\mathcal{K}(\mathcal{H}_{\pi}).
\end{equation}
\ethm
%%--------------------------------------------------------------------
\prf
Observe that it is enough to show that (\ref{eqn:qnzero-17}) holds for the
generating elements $\alpha_{q}$ and $\beta_{q}$ of $A_{q}$. Let us first prove that
%%--------------------------------------------------------------------
\begin{equation}\label{eqn:qnzero-3}
U\lambda_{q}(\alpha_{q})U^{*}-I\otimes \pi_{q}(\alpha_{q})\in \mathscr{T}\otimes\mathcal{K}(\mathcal{H}_{\pi}).
\end{equation}
%%--------------------------------------------------------------------
Note that
%%-------------------------------------
\begin{IEEEeqnarray*}{rCl}
\IEEEeqnarraymulticol{3}{l}{U\lambda_{q}(\alpha_{q})U^{*}e(r,s,t)}\\
 &=& (-1)^{t_{-}}U\lambda_{q}(\alpha_{q})
       e^{\halfof{r+s+|t|}}_{\halfof{-r+s-t},\halfof{-r+s+t}}\\
 &=& (-1)^{t_{-}}U \Bigl(q^{2s+|t|+1}
 \frac{g(r+t_{-}+1)g(r+t_{+}+1)}
 {g(r+s+|t|+1)g(r+s+|t|+2)}
   \;e^{\halfof{r+s+|t|+1}}_{\halfof{-r+s-t-1},\halfof{-r+s+t-1}}\\
  && \qquad + \>
  \frac{g(s+t_{+})g(s+t_{-})}
  {g(r+s+|t|)g(r+s+|t|+1)}
    \;e^{\halfof{r+s+|t|-1}}_{\halfof{-r+s-t-1},\halfof{-r+s+t-1}}\Bigr)
\end{IEEEeqnarray*}
%%-------------------------------------
Since $\max\{\halfof{-r+s-t-1},\halfof{-r+s+t-1}\}=\halfof{-r+s+|t|-1}$
and $\min\{\halfof{-r+s-t-1},\halfof{-r+s+t-1}\}=\halfof{-r+s-|t|-1} $,
we have
%%-------------------------------------
\begin{IEEEeqnarray*}{rCl}
\IEEEeqnarraymulticol{3}{l}{U\lambda_{q}(\alpha_{q})U^{*}e(r,s,t)}\\
  &=& q^{2s+|t|+1}
  \frac{g(r+t_{-}+1)g(r+t_{+}+1)}
  {g(r+s+|t|+1)g(r+s+|t|+2)}
    e(r+1,s,t)\\
  && \qquad\qquad +\> \frac{g(s+t_{+})g(s+t_{-})}
   {g(r+s+|t|)g(r+s+|t|+1)}
  e(r,s-1,t).\yesnumber \label{eqn:qnzero-4}
\end{IEEEeqnarray*}
%%-------------------------------------
For the representation $\pi$, one has
%%-------------------------------------
\begin{equation}\label{eqn:qnzero-5}
(I\otimes\pi_{q}(\alpha_{q}))e(r,s,t)
 = g(s)e(r,s-1,t).
\end{equation}
%%-------------------------------------
From (\ref{eqn:qnzero-4}) and (\ref{eqn:qnzero-5}), we have
\begin{IEEEeqnarray*}{rCl}
%%-------------------------------------
\IEEEeqnarraymulticol{3}{l}{
\Bigl(U\lambda_{q}(\alpha_{q})U^{*}-I\otimes
  \pi_{q}(\alpha_{q})\Bigr) e(r,s,t)}\\
 &=& q^{2s+|t|+1}
  \frac{g(r+t_{-}+1)g(r+t_{+}+1)}
  {g(r+s+|t|+1)g(r+s+|t|+2)}
    e(r+1,s,t)\\
  && \qquad\qquad +\> \Bigl(\frac{g(s+t_{+})g(s+t_{-})}
   {g(r+s+|t|)g(r+s+|t|+1)} -
    g(s)\Bigr)
  e(r,s-1,t).\\
  &=& \bigl(R_{1}(S^{*}\otimes I\otimes I)+R_{2}(I\otimes S \otimes I)\bigr) e(r,s,t),\yesnumber \label{eqn:qnzero-6}
\end{IEEEeqnarray*}
%%-------------------------------------
where $S$ is the left shift on $\ell^{2}(\mathbb{N})$  and
%%-------------------------------------
\begin{IEEEeqnarray}{rCl}
   R_{1} e_{r,s,t} &=& q^{2s+|t|+1}
  \frac{g(r+t_{-}+1)g(r+t_{+}+1)}
  {g(r+s+|t|+1)g(r+s+|t|+2)}
    e(r,s,t)  ,\label{eqn:qnzero-7}\\
   R_{2} e_{r,s,t} &=&  \Bigl(\frac{g(s+t_{+}+1)g(s+t_{-}+1)}
   {g(r+s+|t|+1)g(r+s+|t|+2)} -
    g(s+1)\Bigr) e(r,s,t).\label{eqn:qnzero-8}
\end{IEEEeqnarray}
%%-------------------------------------
Define operators $R_{3}$ and $R_{4}$ on $\mathcal{H}_{\pi}$ as follows:
%%-------------------------------------
\begin{IEEEeqnarray}{rCl}
   R_{3} e_{s,t} &=& q^{2s+|t|+1} e(s,t),\label{eqn:qnzero-9}\\
   R_{4} e_{s,t} &=&  g(s+1)(g(s+|t|+1)-1)e(s,t).\label{eqn:qnzero-10}
\end{IEEEeqnarray}
%%-------------------------------------
It is straightforward to check that $R_{3}$ and $R_{4}$ are compact operators. We
will show that $R_{1}-(I\otimes R_{3})$ and $R_{2}-(I\otimes R_{4})$ are both
compact operators on $\mathcal{H}_{mult}\otimes \mathcal{H}_{\pi}$. It will then
follow from (\ref{eqn:qnzero-6}) that
$U\lambda_{q}(\alpha_{q})U^{*}-I\otimes
  \pi_{q}(\alpha_{q})$ is in $ \mathscr{T}\otimes\mathcal{K}(\mathcal{H}_{\pi})$.

  Observe that
%%-------------------------------------
\begin{IEEEeqnarray*}{rCl}
  \Bigl(R_{1}-(I\otimes R_{3})\Bigr)e(r,s,t)
    &=& q^{2s+|t|+1}\left(\frac{g(r+t_{-}+1)g(r+t_{+}+1)}
  {g(r+s+|t|+1)g(r+s+|t|+2)}-1\right)e(r,s,t).
\end{IEEEeqnarray*}
%%-------------------------------------
We will now use the following estimates for the proof of this:
%%-------------------------------------
\begin{IEEEeqnarray}{lrCll}
	& \left|1-g(k)\right| &<&  q^{2k} &
	   \text{ for all } k\geq 1.\label{eqn:qnzero-1}\\
	\text{There exists $c>0$ such that } &
	   \left|1-g(k)^{-1}\right| &<& c q^{2k} &
	   \text{ for all } k\geq 1.\label{eqn:qnzero-2}
\end{IEEEeqnarray}
%%-------------------------------------
From the above estimates, it
follows that
%%-------------------------------------
\[
\left|q^{2s+|t|+1}\left(\frac{g(r+t_{-}+1)g(r+t_{+}+1)}
  {g(r+s+|t|+1)g(r+s+|t|+2)}-1\right)\right|
     = O(q^{2r+2s+|t|+1}),
\]
%%-------------------------------------
so that
%%-------------------------------------
\begin{equation}\label{eqn:qnzero-20}
	R_{1}-(I\otimes R_{3})\in
	\mathcal{K}(\mathcal{H}_{mult}\otimes\mathcal{H}_{\pi}).
\end{equation}
%%-------------------------------------
Similarly, we have
%%-------------------------------------
\begin{IEEEeqnarray*}{rCl}
\IEEEeqnarraymulticol{3}{l}{
  \Bigl(R_{2}-(I\otimes R_{4})\Bigr)e(r,s,t)}\\
    &=& \Bigl(\frac{g(s+t_{+}+1)g(s+t_{-}+1)}
  {g(r+s+|t|+1)g(r+s+|t|+2)}\\
   && \qquad\qquad - \> g(s+1)+g(s+1)(1-g(s+|t|+1))\Bigr)e(r,s,t)\\
   &=& \Bigl(\frac{g(s+1)g(s+|t|+1)}
  {g(r+s+|t|+1)g(r+s+|t|+2)}-g(s+1)g(s+|t|+1)\Bigr)e(r,s,t).
\end{IEEEeqnarray*}
Using (\ref{eqn:qnzero-1}) and (\ref{eqn:qnzero-2}), we obtain
\begin{IEEEeqnarray*}{rCl}
\IEEEeqnarraymulticol{3}{l}{\Bigl|\frac{g(s+1)g(s+|t|+1)}
  {g(r+s+|t|+1)g(r+s+|t|+2)} - g(s+1)g(s+|t|+1)\Bigr|}\\
     &\leq &  g(s+1)(1-g(s+|t|+1))
	    \Bigl|1-\frac{1}{g(r+s+|t|+1)g(r+s+|t|+2)}\Bigr|\\
 &=& O(q^{2r+2s+2|t|}).
\end{IEEEeqnarray*}
%%-------------------------------------
Hence
%%-------------------------------------
\begin{equation}\label{eqn:qnzero-21}
	R_{2}-(I\otimes R_{4})\in
	\mathcal{K}(\mathcal{H}_{mult}\otimes\mathcal{H}_{\pi}).
\end{equation}
%%-------------------------------------
Combining (\ref{eqn:qnzero-6}), (\ref{eqn:qnzero-20}) and
(\ref{eqn:qnzero-21}), we have (\ref{eqn:qnzero-3}).

Next let us show that
%%--------------------------------------------------------------------
\begin{equation}\label{eqn:qnzero-11}
U\lambda_{q}(\beta_{q})U^{*}-I\otimes \pi_{q}(\beta_{q})\in \mathscr{T}\otimes\mathcal{K}(\mathcal{H}_{\pi}).
\end{equation}
%%--------------------------------------------------------------------
As before, note that
%%-------------------------------------
\begin{IEEEeqnarray*}{rCl}
\IEEEeqnarraymulticol{3}{l}{U\lambda_{q}(\beta_{q})U^{*}e(r,s,t)}\\
 &=& (-1)^{t_{-}}U\lambda_{q}(\beta_{q})
       e^{\halfof{r+s+|t|}}_{\halfof{-r+s-t},\halfof{-r+s+t}}\\
 &=& (-1)^{t_{-}}U \Bigl(-q^{s+t_{+}}
 \frac{g(r+t_{-}+1)g(s+t_{-}+1)}
 {g(r+s+|t|+1)g(r+s+|t|+2)}
   e^{\halfof{r+s+|t|+1}}_{\halfof{-r+s-t+1},\halfof{-r+s+t-1}}\\
  && \qquad + \> q^{s+t_{-}}
  \frac{g(s+t_{+})g(r+t_{+})}
  {g(r+s+|t|)g(r+s+|t|+1)}
    e^{\halfof{r+s+|t|-1}}_{\halfof{-r+s-t+1},\halfof{-r+s+t-1}}\Bigr).
\end{IEEEeqnarray*}
%%-------------------------------------
Since $\max\{\halfof{-r+s-t+1},\halfof{-r+s+t-1}\}=\halfof{-r+s+|t-1|}$,
we have
%%-------------------------------------
\begin{IEEEeqnarray*}{rCl}
\IEEEeqnarraymulticol{3}{l}{U\lambda_{q}(\beta_{q})U^{*}e(r,s,t)}\\
  &=& \begin{cases}
  -q^{s+t}
  \frac{g(r+1)g(s+1)}
  {g(r+s+t+1)g(r+s+t+2)}
    e(r+1,s+1,t-1) &\\
  \qquad\qquad +\> q^{s} \frac{g(s+t)g(r+t)}
   {g(r+s+t)g(r+s+t+1)}
  e(r,s,t-1) &\text{if }t \geq 1,\\
    &\yesnumber \label{eqn:qnzero-12}\\
  q^{s}
  \frac{g(r+|t|+1)g(s+|t|+1)}
  {g(r+s+|t|+1)g(r+s+|t|+2)}
    e(r,s,t-1) &\\
  \qquad\qquad -\> q^{s+|t|} \frac{g(s)g(r)}
   {g(r+s+|t|)g(r+s+|t|+1)}
  e(r-1,s-1,t-1) &\text{if }t < 1.
  \end{cases}
\end{IEEEeqnarray*}
%%-------------------------------------
For the representation $\pi$,
%%-------------------------------------
\begin{equation}\label{eqn:qnzero-13}
(I\otimes\pi_{q}(\beta_{q}))e(r,s,t)
 = q^{s}e(r,s,t-1).
\end{equation}
%%-------------------------------------
From (\ref{eqn:qnzero-12}) and (\ref{eqn:qnzero-13}), we have
%%-------------------------------------
\begin{IEEEeqnarray*}{rCl}
\IEEEeqnarraymulticol{3}{l}{\Bigl(U\lambda_{q}(\beta_{q})U^{*}-
       (I\otimes\pi_{q}(\beta_{q}))\Bigr)e(r,s,t)}\\
  &=& \begin{cases}
  -q^{s+t}
  \frac{g(r+1)g(s+1)}
  {g(r+s+t+1)g(r+s+t+2)}
    e(r+1,s+1,t-1) &\\
  \qquad\qquad +\> q^{s} \Bigl(\frac{g(s+t)g(r+t)}
   {g(r+s+t)g(r+s+t+1)}-1\Bigr)
  e(r,s,t-1) &\text{if }t \geq 1,\\
    &\\
  q^{s}
  \Bigl(\frac{g(r+|t|+1)g(s+|t|+1)}
  {g(r+s+|t|+1)g(r+s+|t|+2)}-1\Bigr)
    e(r,s,t-1) &\\
  \qquad\qquad -\> q^{s+|t|} \frac{g(s)g(r)}
   {g(r+s+|t|)g(r+s+|t|+1)}
  e(r-1,s-1,t-1) &\text{if }t < 1
  \end{cases}\\
  &=&
  (T_{1}(S^{*}\otimes S^{*}\otimes S + S\otimes S\otimes S) + T_{2}(I\otimes I\otimes S))e(r,s,t),\yesnumber\label{eqn:qnzero-14}
\end{IEEEeqnarray*}
%%-------------------------------------
where
$S$ denotes the left shift and
%%-------------------------------------
\begin{IEEEeqnarray}{rCl}
   T_{1} e_{r,s,t}
    &=& \begin{cases}
	     -q^{s+|t|}
  \frac{g(r)g(s)}{g(r+s+|t|)g(r+s+|t|+1)}
    e(r,s,t)  &\text{if } (r,s)\neq (0,0), t\geq 0,\\
  	     -q^{s+|t|}
    \frac{g(r+1)g(s+1)}{g(r+s+|t|+1)g(r+s+|t|+2)}
      e(r,s,t)  &\text{if } (r,s)\neq (0,0), t< 0,\\
	  0  &\text{if } (r,s)= (0,0)
	  \end{cases}\label{eqn:qnzero-15}\\
   T_{2} e_{r,s,t}
     &=&  \begin{cases}
	    q^{s}\Bigl(
		\frac{g(r+|t|+1)g(s+|t|+1)}{g(r+s+|t|+1)g(r+s+|t|+2)}-1\Bigr)
		e(r,s,t) &\text{if }  t\geq 0,\\
    q^{s}\Bigl(
	\frac{g(r+|t|)g(s+|t|)}{g(r+s+|t|)g(r+s+|t|+1)}-1\Bigr)
	e(r,s,t) &\text{if }  t< 0.\label{eqn:qnzero-16}
	       \end{cases}
\end{IEEEeqnarray}
%%--------------------------------------------------------------------
Define operators $T_{3}$ and $T_{4}$ on $\mathcal{H}_{\pi}$ as follows:
%%-------------------------------------
\begin{IEEEeqnarray}{rCl}
   T_{3} e_{s,t} &=& \begin{cases}
	     -q^{s+|t|}g(s)
    e(r,s,t)  &\text{if } (r,s)\neq (0,0), t\geq 0,\\
  	     -q^{s+|t|}g(s+1)
      e(r,s,t)  &\text{if } (r,s)\neq (0,0), t< 0,\\
	  0  &\text{if } (r,s)= (0,0),
	  \end{cases}\\
   T_{4} e_{s,t} &=&  \begin{cases}
	    q^{s}\Bigl(g(s+|t|+1)-1\Bigr)
		e(r,s,t) &\text{if }  t\geq 0,\\
    q^{s}\Bigl(g(s+|t|)-1\Bigr)
	e(r,s,t) &\text{if }  t< 0.
	       \end{cases}
\end{IEEEeqnarray}
%%-------------------------------------
It is straightforward to check that $T_{3}$ and $T_{4}$ are compact operators on
$\mathcal{H}_{\pi}$. We will show that $T_{1}-(I\otimes T_{3})$ and
$T_{2}-(I\otimes T_{4})$ are both compact operators on $\mathcal{H}_{mult}\otimes
\mathcal{H}_{\pi}$. It will then follow from (\ref{eqn:qnzero-14}) that
$U\lambda_{q}(\beta_{q})U^{*}-I\otimes
  \pi_{q}(\beta_{q})$ is in $ \mathscr{T}\otimes\mathcal{K}(\mathcal{H}_{\pi})$.

  Observe that
%%-------------------------------------
\begin{IEEEeqnarray*}{rCl}
\IEEEeqnarraymulticol{3}{l}{
  \Bigl(T_{1}-(I\otimes T_{3})\Bigr)e(r,s,t)}\\
    &=& \begin{cases}
	     -q^{s+|t|}g(s)
   \Bigl(\frac{g(r)}{g(r+s+|t|)g(r+s+|t|+1)}-1\Bigr)
    e(r,s,t)  &\text{if } (r,s)\neq (0,0), t\geq 0,\\
  	     -q^{s+|t|}g(s+1)
    \Bigl(\frac{g(r+1)}{g(r+s+|t|+1)g(r+s+|t|+2)}-1\Bigr)
      e(r,s,t)  &\text{if } (r,s)\neq (0,0), t< 0,\\
	  0  &\text{if } (r,s)= (0,0).
	  \end{cases}
\end{IEEEeqnarray*}
%%-------------------------------------
From the estimates in (\ref{eqn:qnzero-1}) and (\ref{eqn:qnzero-2}), it
follows that the right hand side above is $O(q^{r+s+|t|})$.
Hence
%%-------------------------------------
\begin{equation}\label{eqn:qnzero-18}
	T_{1}-(I\otimes T_{3})\in
	\mathcal{K}(\mathcal{H}_{mult}\otimes\mathcal{H}_{\pi}).
\end{equation}
%%-------------------------------------

Similarly, one has
%%-------------------------------------
\begin{IEEEeqnarray*}{rCl}
\IEEEeqnarraymulticol{3}{l}{
  \Bigl(T_{2}-(I\otimes T_{4})\Bigr)e(r,s,t)}\\
    &=& \begin{cases}
	    q^{s}g(s+|t|+1)\Bigl(
		\frac{g(r+|t|+1)}{g(r+s+|t|+1)g(r+s+|t|+2)}-1\Bigr)
		e(r,s,t) &\text{if }  t\geq 0,\\
    q^{s}g(s+|t|)\Bigl(
	\frac{g(r+|t|)}{g(r+s+|t|)g(r+s+|t|+1)}-1\Bigr)
	e(r,s,t) &\text{if }  t< 0.
	       \end{cases}
\end{IEEEeqnarray*}
%%-------------------------------------
Using (\ref{eqn:qnzero-1}) and (\ref{eqn:qnzero-2}), we conclude that the
right hand side is $O(q^{r+s+|t|})$. Therefore
%%-------------------------------------
\begin{equation}\label{eqn:qnzero-19}
	T_{2}-(I\otimes T_{4})\in
	 \mathcal{K}(\mathcal{H}_{mult}\otimes\mathcal{H}_{\pi}).
\end{equation}
%%-------------------------------------
Thus combining (\ref{eqn:qnzero-14}), (\ref{eqn:qnzero-18}) and
(\ref{eqn:qnzero-19}), we get (\ref{eqn:qnzero-11}).
\qed

%%--------------------------------------------------------------------
\section{Applications}
%%--------------------------------------------------------------------
We now describe three different contexts in noncommutative topology and geometry
where the main result of the last section can be used.
%%--------------------------------------------------------------------
\subsection{$KK$-groups}
%%--------------------------------------------------------------------
The first is in the context of $KK$ theory. The approximate
equivalence tells us that
%%-------------------------------------
\begin{equation}\label{eq:con-rem1}
U\lambda_{q}(a)U^{*}-I\otimes \pi_{q}(a)\in
\mathscr{T}\otimes\mathcal{K}(\mathcal{H}_{\pi}) \text{ for all }a\in A_{q}.
\end{equation}
%%-------------------------------------
Since $I\otimes
\pi_{q}(a) \in \mathscr{T}\otimes \mathcal{L}(\mathcal{H}_{\pi})$ and
$\mathscr{T}\otimes \mathcal{L}(\mathcal{H}_{\pi})\subseteq
M(\mathscr{T}\otimes\mathcal{K}(\mathcal{H}_{\pi}))$, one gets
%%-------------------------------------
\begin{equation}\label{eq:con-rem2}
U\lambda_{q}(a)U^{*}\in M(\mathscr{T}\otimes\mathcal{K}(\mathcal{H}_{\pi}))
   \text{ for all }a\in A_{q}.
\end{equation}
%%-------------------------------------
It now follows that the pair $(U\lambda_{q}(\cdot)U^{*}, I\otimes\pi_{q}(\cdot))$
gives a quasihomomorphism in the Cuntz description (see \cite{Cun-1987np}) of the
$KK$ group $KK(A_{q},\mathscr{T})$. In other words, we have the following theorem.

%%-------------------------------------
\bthm
There is a unitary $U:L^{2}(SU_q(2))\to
\mathcal{H}_{mult}\otimes\mathcal{H}_{\pi}$ such that %the pair
$(U\lambda_{q}(\cdot)U^{*}, \pi_{q})$ is a quasihomomorphism and gives a $KK$-class
in the group $KK(A_{q},\mathscr{T})$.
\ethm
%%-------------------------------------

%%--------------------------------------------------------------------
\subsection{Fredholm representation of $\widehat{SU_{q}(2)}$}
%%--------------------------------------------------------------------
Next we will show that the above approximate equivalence also gives an example of a
Fredholm representation of the dual quantum group $\widehat{SU_{q}(2)}$. Let us
first recall the notion of a Fredholm representation of a discrete group.
%%-------------------------------------
\bdfn (Mishchenko \cite{Mis-1975pz})
Let $\Gamma$ be a discrete group. A pair of unitary representations $\pi_{1}$ and
$\pi_{2}$ of $\Gamma$ acting on a Hilbert space $\mathcal{H}$ together with a
Fredholm operator $F\in\mathcal{L}(\mathcal{H})$ is called a Fredholm
representation of $\Gamma$ if $F\pi_{1}(g)-\pi_{2}(g)F\in \mathcal{K}(\mathcal{H})$
for all $g\in\Gamma$.
\edfn
%%-------------------------------------
This notion admits an immediate extension to discrete quantum groups where the
representations $\pi_{1}$ and $\pi_{2}$ will be elements of $M(\hat{A}\otimes
\mathcal{K}(\mathcal{H}))$, where $\hat{A}$ stands for the $C^{*}$-algebra
associated with the quantum group and the condition $F\pi_{1}(g)-\pi_{2}(g)F\in
\mathcal{K}(\mathcal{H})$ gets replaced by $(I\otimes F)\pi_{1}-\pi_{2}(I\otimes F)
\in \hat{A}\otimes\mathcal{K}(\mathcal{H})$.

We now extend this notion further and incorporate coefficients
from a $C^*$-algebra $B$.
%%-------------------------------------
\bdfn
Assume that the Hilbert space $\mathcal{H}$ on which the two representations act is
of the form $\mathcal{H}_{1}\otimes \mathcal{H}_{2}$, and that $B$ is a
$C^{*}$-subalgebra of $\mathcal{L}(\mathcal{H}_{1})$ such that one has
\[
\pi_{1}, \pi_{2} \in M(\hat{A}\otimes B\otimes \mathcal{K}(\mathcal{H}_{2})).
\]
We call $(\pi_{1},\pi_{2},F)$ a \textbf{Fredholm representation with coefficients
in $B$} if
\[
(I\otimes
F)\pi_{1}-\pi_{2}(I\otimes F) \in \hat{A}\otimes\mathcal{K}(B\otimes\mathcal{H}),
\]
where
$B\otimes\mathcal{H}$ denotes the standard Hilbert $B$-module $\ell^{2}(B)$.
\edfn
%%-------------------------------------

Given a representation $\pi$ of the $C^{*}$-algebra $A\equiv C(G)$ for a compact
quantum group $G$, one can associate a co-representation of the dual quantum group
$\widehat{G}$ with its opposite coproduct as follows. Let $\Gamma$ be the set of
equivalence classes of unitary representations of the compact quantum group $G$.
Denote by $u^{(\gamma)}$ the irreducible unitary representation indexed by
$\gamma\in\Gamma$ acting on a finite dimensional Hilbert space
$\mathcal{H}_{\gamma}$. Thus $u^{(\gamma)}\in
\mathcal{L}(\mathcal{H}_{\gamma})\otimes A$. Let $\hat{A}$ denote the
$c_{0}$-direct sum $\oplus_{\gamma\in\Gamma}\mathcal{L}(\mathcal{H}_{\gamma})$.
Then $u:=\oplus_{\gamma}u^{(\gamma)}$ is a unitary element of $M(\hat{A}\otimes A)$
and satisfies the following two identities:
%%-------------------------------------
\[
(id\otimes \Delta_{G})u=u_{1,2}u_{1,3},\quad
  (\Delta_{\hat{G}}\otimes id)u= u_{2,3}u_{1,3}.
\]
%%-------------------------------------
Now define $w_{\pi}:=(id\otimes\pi)u$. Then $w_{\pi}$ is a unitary element of
$M(\hat{A}\otimes \mathcal{H})$ and satisfies $(\Delta_{\hat{G}}\otimes id)w_{\pi}=
(w_{\pi})_{2,3}(w_{\pi})_{1,3}$. Thus the element $w_{\pi}$ gives a unitary
(co-)representation of the dual quantum group $\hat{G}$, which has the same set of
intertwiners as $\pi$.

Next, note that
%%-------------------------------------
\[
\mathscr{T}\otimes\mathcal{K}(\mathcal{H}_{\pi})\cong
\mathcal{K}(\mathscr{T}\otimes \mathcal{H}_{\pi}), \qquad
M(\mathscr{T}\otimes\mathcal{K}(\mathcal{H}_{\pi}))\cong
\mathcal{L}(\mathscr{T}\otimes \mathcal{H}_{\pi}),
\]
%%-------------------------------------
where $\mathcal{K}(\mathscr{T}\otimes \mathcal{H}_{\pi})$ and
$\mathcal{L}(\mathscr{T}\otimes \mathcal{H}_{\pi})$ denote the space of compact
operators and the space of adjointable operators respectively on the Hilbert
$\mathscr{T}$-module $\mathscr{T}\otimes \mathcal{H}_{\pi}$. Thus we have for all
$a\in A_{q}$,
%%-------------------------------------
\begin{equation}\label{eq:con-rem3}
U\lambda_{q}(a)U^{*},\; I\otimes \pi_{q}(a)
   \in \mathcal{L}(\mathscr{T}\otimes \mathcal{H}_{\pi}),\qquad
 U\lambda_{q}(a)U^{*}- I\otimes \pi_{q}(a)
   \in \mathcal{K}(\mathscr{T}\otimes \mathcal{H}_{\pi}).
\end{equation}
%%-------------------------------------
If we denote by $w_{\lambda}$ and $w_{\pi}$ the co-representations of the dual
quantum group $\widehat{SU_{q}(2)}$ corresponding to the representations
$U\lambda_{q}(\cdot)U^{*}$ and $I\otimes \pi_{q}(\cdot)$ respectively of the
$C^{*}$-algebra $A_{q}$, and by $\widehat{A_{q}}$ the $C^{*}$-algebra associated
with the dual $\widehat{SU_{q}(2)}$, then
$w_{\lambda}$ and $w_{\pi}$ are unitary elements of $M(\widehat{A_{q}}\otimes
\mathcal{K}(\ell^{2}(\mathbb{N})\otimes \mathcal{H}_{\pi}))$. It follows from
(\ref{eq:con-rem1}) that
%%-------------------------------------
\[
 w_{\lambda}- w_{\pi}\in \widehat{A_{q}}\otimes
   \mathcal{K}(\mathscr{T}\otimes \mathcal{H}_{\pi}).
\]
%%-------------------------------------
Thus $(w_{\lambda}, w_{\pi}, I)$ gives a Fredholm representation of the dual $\widehat{SU_{q}(2)}$ with coefficients in the Toeplitz algebra.

%%--------------------------------------------------------------------
\subsection{Spectral triples}
%%--------------------------------------------------------------------
The spectral triples corresponding to the two representations $\lambda_{q}$ and
$\pi_{q}$ given in \cite{ChaPal-2003ab} and \cite{ChaPal-2003aa} are
$\left(\ell^{2}(\mathbb{N} \times \mathbb{Z}), \pi_{q}, D_{\pi}\right)$ and
$\left(L^{2}\left(SU_{q}(2)\right), \lambda_{q}, D_{\lambda}\right)$, and one has
%%-------------------------------------
\begin{IEEEeqnarray*}{rCl}
	U \lambda_{q}(a) U^{*} &=& I \otimes \pi_{q}(a) \text{ modulo }
	 \mathscr{T} \otimes \mathcal{K}(\ell^{2}(\mathbb{N} \times \mathbb{Z})),\\
	U\left|D_{\lambda}\right| U^{*} &=& I \otimes\left|D_{\pi}\right|+N \otimes I,\\
	U(\sgn D_{\lambda}) U^{*} &=& P_{0} \otimes(\sgn D_{\pi}) - (I-P_{0}) \otimes I.
\end{IEEEeqnarray*}
%%-------------------------------------
One can now establish that the differences $U \lambda_{q}(a) U^{*} - I \otimes
\pi_{q}(a)$ belong to certain finer ideals of the compacts for specific collections
of elements $a\in A_{q}$ and this enables one to use the knowledge of regularity
and dimension spectrum of the spectral triple in \cite{ChaPal-2003ab} to draw
similar conclusions on the spectral triple in \cite{ChaPal-2003aa}.

\brmrk
For odd dimensional quantum spheres $S_{q}^{2n+1}$, which are the homogeneous
spaces $SU_{q}(n+1)/SU_{q}(n)$, using the  $C^{*}$-algebra $C(S_{q}^{2n+1})$  at
$q=0$ one can prove an approximate
equivalence between the representation of $C(S_{q}^{2n+1})$ on the $L^{2}$-space
of the invariant state and a faithful representation that is easier to work with
computationally. In this case, if one works with the basis for the $L^{2}$-space
that comes from the GT basis for the $L^{2}$-space of the Haar State on
$SU_{q}(n+1)$, then the resulting $\Gamma$, though more complicated than the
$SU_{q}(2)$ case, is very similar and can be handled in essentially the same way to
arrive at the unitary. The approximate equivalence can then be derived  in a
very similar manner as in the present case.
\ermrk

\brmrk
It is worthwhile to make a remark at this point about crystallisations of the
$C^*$-algebras $C(G_{q})$ where $G_{q}$ stands for the $q$-deformation of a
connected simply connected compact Lie group. In two recent papers, Pal \&
Giri \cite{GirPal-2022tv} and Matassa \& Yuncken \cite
{MatYun-2022aa} introduced the notion of crystallisation of these
$C^*$-algebras. The paper \cite{GirPal-2022tv} then focuses on studying the
irreducible representations of the crystallised $C^*$-algebra in the type $A$
case, while in \cite{MatYun-2022aa} the authors exploit the link with
Kashiwara's crystal basis theory to study certain properties of the
crystallised $C^*$-algebras. Our present paper is an illustration of how these
crystallised $C^*$-algebras can be useful in the type $A_{1}$ case. One can
hope to obtain similar decomposition results for  higher rank cases and for
deformations of Lie groups of other types ($B$, $C$, $D$ etc.) using the
crystallised $C^*$-algebras introduced in \cite{GirPal-2022tv} and  \cite
{MatYun-2022aa}, which in turn should pave the way for a detailed study of the
Neshveyev-Tuset Dirac operators for $q$-deformations of Lie groups of higher
ranks.
\ermrk

%%--------------------------------------------------------------------
\bibliographystyle{plain}
% \bibliography{libminimal}
%%--------------------------------------------------------------------
%%--------------------------------------------------------------------
%%--------------------------------------------------------------------

%%--------------------------------------------------------------------
% \hrulefill
%
% \noindent{\sc Partha Sarathi Chakraborty}
% (\texttt{parthacsarathi.isi.smu@gmail.com})\\
% {\footnotesize Indian Statistical
% Institute,  Kolkata, INDIA}\\[1ex]
% {\sc Arup Kumar Pal} (\texttt{arupkpal@gmail.com})\\
%          {\footnotesize Indian Statistical
% Institute, Delhi Centre, INDIA}
%%--------------------------------------------------------------------
\end{document}